\documentclass[a4paper,12pt,leqno]{article}
\usepackage{graphicx}
\usepackage{amsmath}
\usepackage{amssymb}
\usepackage{amstext}
\usepackage{amsxtra}
\usepackage[mathscr]{eucal}

\makeatletter
\@addtoreset{equation}{section}
\renewcommand\theequation
  {\ifnum \c@section>\z@ \thesection.\fi \@arabic\c@equation}
\makeatother

\newcommand{\bR}{{\mathbb{R}}}

\newcommand{\bN}{{\mathbb{N}}}

\newcommand{\cK}{{\cal K}}

\newcommand{\sG}{{\mathscr{G}}}

\newcommand{\sI}{{\mathscr{I}}}
\newcommand{\sJ}{{\mathscr{J}}}
\newcommand{\sK}{{\mathscr{K}}}
\newcommand{\sL}{{\mathscr{L}}}
\newcommand{\sM}{{\mathscr{M}}}
\newcommand{\sP}{{\mathscr{P}}}

\newcommand{\II}{I\kern -1mm I}
\newcommand{\III}{I\kern -1mm I\kern -1mm I}
\newcommand{\al}{\alpha}
\newcommand{\be}{\beta}
\newcommand{\e}{\varepsilon}
\newcommand{\ga}{\gamma}

\newcommand{\p}{\partial}

\newcommand{\om}{\Omega}
\newcommand{\oom}{\overline\Omega}
\newcommand{\pom}{\partial\Omega}

\newcommand{\whvphi}{\widehat\varphi}

\newcommand{\os}{\overline s}

\newcommand{\ot}{\overline t}
\newcommand{\ou}{\overline u}
\newcommand{\uu}{\underline u}
\newcommand{\wu}{\widetilde u}
\newcommand{\oU}{\overline U}
\newcommand{\uU}{\underline U}

\newcommand{\uv}{\underline v}

\newcommand{\oW}{\overline W}

\newcommand{\ox}{\overline x}
\newcommand{\ux}{\underline x}
\newcommand{\oy}{\overline y}

\newcommand{\ds}{\displaystyle}
\newcommand{\pprime}{{\prime\prime}}

\newcommand{\sh}{\kern 0.5mm{\rm sh}}
\newcommand{\ch}{{\rm ch}}
\newtheorem{theorem}{Theorem}[section]
\newtheorem{remark}{Remark}[section]
\newtheorem{proposition}{Proposition}[section]

\newtheorem{lemma}{Lemma}[section]

\title{\bf Convergence of the Approximation scheme to American option 
pricing via the discrete Morse semiflow}
\date{}
\author{
\medskip
Katsuyuki Ishii \\
Graduate School of Maritime Sciences \\
Kobe University \\
\medskip
Higashinada, Kobe 658-0022, JAPAN \\
e-mail:ishii@maritime.kobe-u.ac.jp \\
\medskip
and \\
\medskip
Seiro Omata \\
Graduate School of Natural Science and Technology \\
Kanazawa University \\
Kanazawa 920-1192, JAPAN \\
e-mail:omata@kenroku.kanazawa-u.ac.jp
}

\begin{document}

\maketitle

\begin{quote}
{{\bfseries Abstract.}
We consider the approximation scheme of the American call option 
via the discrete Morse semiflow.  It is the minimizing scheme of 
a time-semidiscretized variational functional.  In this paper we 
obtain a rate of convergence of approximate solutions.  
In addition, the convergence of approximate free boundaries 
is proved.}
\end{quote}

\section{Introduction}

In this paper we consider an approximation scheme to the following 
obstacle problem: 
\begin{equation}
\label{bs}
   \left\{\begin{array}{ll}
   \ds{\min\left\{-C_\tau -\frac{\sigma^2}{2}S^2C_{SS}-(r-q)SC_S
  +rC,C-\Phi\right\}=0} & {\rm in\ }(0,T)\times(0,+\infty),\\
   C(T,S)=\Phi(S):=\max(S-K,0) & {\rm for}\ S\in[0,+\infty),\\
   C(\tau,0)=0 & {\rm for\ }\tau\in(0,T),\\
   \end{array}\right.
\end{equation}
The above equation is called the Black-Scholes equation for the 
American call option.  Here $ C=C(\tau,S)$ is the option price, the 
positive constants $ \sigma $, $ r $, $ q $, $ K $ denote, respectively, 
the volatility, the interest rate, the dividend and the strike price.  
Throughout this paper we set $ \sigma=\sqrt{2} $ for simplicity and 
assume $ r>q $.  

It is known that (\ref{bs}) has a unique free boundary.  Indeed, 
by van Moerbeke \cite{vm;76}, we see that under the assumption $ r>q $, 
there exists a unique $ S^*\in C^1(0,T) $ such that 
$ (S^*)^\prime\leq 0 $ and 
\begin{eqnarray}
\nonumber
  & & \{S\;|\;C(\tau,S)>\Phi(S)\}=(-\infty,S^*(\tau)),\ 
      \{S\;|\;C(\tau,S)=\Phi(S)\}=[S^*(\tau),+\infty) \\
\nonumber
  & & \hspace{100mm}  {\rm for\ each\ }\tau\in(0,T),\\
\label{eq101}
   & & \lim_{\tau\nearrow T}S^*(\tau)=\frac{rK}{q}\;(>K).  
\end{eqnarray}
See Wilmott - Howison - Dewynne \cite[p.124]{wi;ho;de;95} for the formal 
derivation of (\ref{eq101}).  The family $ \{S^*(\tau)\}_{0<\tau<T} $ 
is called the free boundary or the optimal exercise boundary; for each 
$ \tau\in (0,T) $, $ S^*(\tau) $ indicates the value of the current 
stock price under which the holder of the option should (optimally) 
exercise it.  

From the viewpoint of mathematical finance, it would be very convenient 
to obtain $ \{S^*(\tau)\}_{0<\tau<T} $ explicitely.  However, this 
seems to be very difficult and thus many people have studied numerical 
schemes for (\ref{bs}), especially to approximate the free boundary.  
Brennan - Schwartz \cite{bre;sch;77} introduced a fully implicit 
difference scheme for the American put option and obtained a numerical 
solution.  The convergence of their scheme was proved by 
Jaillet - Lamberton - Lapeyre \cite{jai;lam;lap;90} in the framework 
of variational inequalities.  Lamberton \cite{lam;93} considered the 
binomial tree method and the finite difference one to approximate 
$ \{S^*(\tau)\}_{0<\tau<T} $ and showed the convergence of the 
approximate free boundary by the probabilistic argument and the 
analytical one.  He also obtained in \cite{lam;98, lam;02} some error 
estimates for the stochastic approximation to the optimal stopping 
problems including the American options.  Amin - Khanna 
\cite{ami;kha;94} treated a discrete time model for the American 
option and proved the convergence of the discrete American option 
value to the continuous one.  Jiang - Dai \cite{jia;dai;04} obtained 
similar results to those in \cite{lam;93} by the method of 
viscosity solutions.  Omata - Iwasaki - Nakane - Xiong - Sakuma 
\cite{om;iw;na;xi;sak;03} proposed an approximation scheme to 
(\ref{bs}) different from the above ones and obtained a numerical result.  

The approximation scheme by \cite{om;iw;na;xi;sak;03} is based on the 
discrete Morse semiflow (DMS), consisting of the minimization of 
a time-semidiscretized variational functional.  The DMS was first used 
by Rektorys \cite{re;71} to obtain the solutions of linear parabolic 
equations.  Kikuchi \cite{ki;91, ki;94} applied the DMS to construct 
the solutions of parabolic equations associated with a variational 
functional of a harmonic map type.  Besides, in 
\cite{ta;94, na;ta;95, na;ta;96, na;ta;96-2} Nagasawa and Tachikawa 
used the DMS to show the existence and asymptotic behavior of solutions 
of some semilinear hyperbolic systems.  
Nagasawa - Omata \cite{na;om;93} considered the behavior of the DMS 
for a free boundary problem.  Some applications of the DMS to numerical 
analysis have been treated in Omata \cite{om;96, om;97}, 
Omata - Okamura - Nakane \cite{om;ok;na;99} 
and Omata - Iwasaki - Kawagoe \cite{om;iw;ka;01}.  

The purpose of this paper is to discuss the convergence of the 
approximation scheme by \cite{om;iw;na;xi;sak;03}.  Our results are 
a rate of convergence of the approximate solutions and the convergence 
of the approximate free boundary.  The former result is obtained by 
applying the rate of convergence of product formula for semigroups 
by Bentkus - Paulauskas \cite{ben;pau;04} and the precise comparison 
argument for viscosity solutions by Ishii - Koike \cite{is;ko;91-2}, 
in which they obtained the rate of convergence in elliptic singular 
perturbations.  The latter one is proved by the limit operation of 
viscosity solutions due to Barles - Perthame \cite{ba;pe;87, ba;pe;88}.  

This paper is organized as follows.  In Section 2 we introduce the 
approximation scheme by \cite{om;iw;na;xi;sak;03} and state the main 
results.  In Section 3 we discuss the solutions of (\ref{bs04}) below.  
In Subsection 3.1 we briefly prove the existence and uniqueness of 
solutions.  In Subsection 3.2 we show some properties of solutions.  
Section 4 is devoted to the DMS associated with (\ref{bs04}) (DMS-BS 
for short) and the free boundary of the DMS-BS.  In Subsection 4.1 
we derive some estimates for the DMS-BS.  Subsection 4.2 is devoted 
to the existence and uniqueness of the free boundary of the DMS-BS.  
In Subsection 4.3, we give a proof of Theorem \ref{th406} in 
subsection 4.1, an estimate of the difference of the DMS-BS.  
In Section 5 we prove our main results.  Section 6 is the Appendix.  
In Subsection 6.1 we discuss the formal asymptotic expansion of 
an ODE related to (\ref{ellvi02}) below.  This expansion is used to 
construct sub- and supersolutions of (\ref{ellvi}) and (\ref{ellvi02}).  
In Subsection 6.2 we give an estimate for some coefficients appearing 
in the estimate of Theorem \ref{th407}.  

In the following of this paper, we denote by $ C $ various constants 
depending only on known ones.  The value of $ C $ may vary from line 
to line.  

\section{Approximation scheme and Main Results}

In this section we state the approximation scheme by 
\cite{om;iw;na;xi;sak;03} and our main results.  

We reformulate (\ref{bs}) in the following way.   Put $ t:=T-\tau $, 
$ x:=\log(S/K) $, $ \al:=(r-q-1)/2 $ and 
$ U(t,x):=S^\al C(\tau,S)/K^{\al+1} $.  Then (\ref{bs}) turns to 
\begin{equation}
\label{bs03}
   \left\{\begin{array}{ll}
   \min\left\{U_t-U_{xx}
  +\be U,U-\varphi\right\}=0 &{\rm in\ }(0,T)\times\bR,\\
   U(0,x)=\varphi(x):=e^{\al x}\max(e^x-1,0) & {\rm for}\ x\in\bR,\\
  U(t,x)\longrightarrow 0 \quad (x\to -\infty)
      & {\rm for\ }t\in(0,T),\\ 
   \end{array}\right.
\end{equation}
where $ \be:=\al^2+r $.  From the viewpoint of the numerical analysis, 
we had better restrict the problem (\ref{bs03}) on a bounded interval 
with respect to $ x $.  This restriction seems to be reasonable.  
Because the free boundary $ \{\log(S^*(T-t)/K)\}_{0<t<T} $ for 
(\ref{bs03}) is bounded for each $ T>0 $ and it is easily seen that 
$$
   U(t,x)=O(e^{\ga x})\quad{\rm as\ }x\to -\infty\ 
   \quad{\rm for\ all\ }t\in (0,T)\ {\rm and\ some\ }\ga>0.  
$$
Hence, putting $ \om:=(-1,1) $, we consider the following problem 
instead of (\ref{bs03}): 
\begin{equation}
\label{bs04}
   \left\{\begin{array}{ll}
   \ds{\min\left\{u_t-u_{xx}
  +\be u,u-\varphi\right\}=0} & {\rm in\ }(0,T)\times\om,\\
   u(0,x)=\varphi(x) & {\rm for}\ 
    x\in\overline\om,\\
  u(t,\pm 1)=\varphi(\pm 1) & {\rm for\ }t\in(0,T).\\ 
   \end{array}\right.
\end{equation}

We assume $ q<r<qe $ and denote by $ \{x^*(t)\}_{0<t<T} $ the free 
boundary for (\ref{bs04}).  Note by \cite{vm;76} that 
\begin{equation}
\label{free}
  x^*\in C^1(0,T),\ (x^*)^\prime(t)\geq 0,\ 
   \lim_{t\searrow 0}x^*(t)=x_0
   :=\log\left(\frac{r}{q}\right)(\in(0,1)).  
\end{equation}

The approximation scheme by \cite{om;iw;na;xi;sak;03} is stated 
as follows.  Fix a time step $ h>0 $.  Put $ u_0:=\varphi $ and let 
$ [r] $ be the Gauss symbol for $ r\in\bR$.  For 
$ m=1,2,\ldots, [T/h] $, we consider the minimization problem of the 
following functional: 
\begin{eqnarray*}
  & & J_m(u):=\frac{1}{2}\int_{\om}
   \left\{\frac{|u-u_{m-1}|^2}{h}+(u_x)^2+\be u^2\right\}dx\quad
    {\rm for\ }u\in\cK,\\
  & & \cK:=\{v\in H^1(\om)\;|\;v-\varphi\in H_0^1(\om),\ 
    v\geq\varphi\ {\rm a.e.\ in\ }\om\}.  
\end{eqnarray*}
We observe by the direct method of calculus of variation that there is 
a unique minimizer $ u_m\in \cK $ of $ J_m $.  Moreover, $ u_m $ 
satisfies the elliptic variational inequality:
\begin{equation}
\label{ellvi}
    \min
   \left\{\frac{u_m-u_{m-1}}{h}-u_{m,xx}+\be u_m,u_m-\varphi\right\}=0
    \quad {\rm in\ }\om, \ 
     u_m(\pm 1)=\varphi(\pm 1).  
\end{equation}
We call the sequence $ \{u_m\}_{m=0}^{ [T/h]} $ the DMS-BS. 
In addition, there is a unique free boundary $ \{x_m\}_{m=0}^{[T/h]} $ 
to the DMS-BS, as will be shown in Subsection 4.2 below.  

Under these settings, we define $ u^h(t,x) $ and $ x^h(t) $ by 
\begin{eqnarray}
\label{approx01}
  & & u^h(t,x):=u_m(x),\ x^h(t):=x_m\\
\nonumber  
  & & \hspace{30mm} {\rm for\ }t\in [mh,(m+1)h),\ x\in\oom\ 
      {\rm and\ }m=0,1,\ldots,[T/h]. 
\end{eqnarray}
Then our main results are stated as follows.  
\begin{theorem}
\label{th201}
Assume $ q<r<qe $.  Then for any $ \delta>0 $, there exist $ K>0 $ and 
$ h_0>0 $ such that for all $ h\in (0,h_0) $, 
$$
  \sup_{(t,x)\in[0,T-\delta]\times\oom}|u^h(t,x)-u(t,x)|
    \leq K\sqrt{h}|\log h|.
$$
\end{theorem}
\begin{theorem}
\label{th202}
Assume $ q<r<qe $ and $ \overline{\{x^*(t)\}_{0<t<T}}\subset\om $.  
Then for any $ \delta>0 $, we have 
$$
   \lim_{h\to 0}\sup_{t\in [0,T-\delta]}|x^h(t)-x^*(t)|=0.  
$$
\end{theorem}
\begin{remark}{\rm 
(1) The assumptions $ r<qe $ and 
$ \overline{\{x^*(t)\}_{0<t<T}}\subset\om $ are technical ones.  Since 
$ \overline{\{x^*(t)\}_{0<t<T}} $ is bounded for each $ T>0 $ 
(cf. \cite{vm;76}), replacing $ \om $ with a larger interval such that 
$ \overline{\{x^*(t)\}_{0<t<T}}\subset\om $, we can show Theorems 
\ref{th201} and \ref{th202}, assuming only $ r>q $.  

(2) The $ |\log h| $ term appear in Theorem \ref{th201} by some 
technical reasons.  
}
\end{remark}

\section{Solutions of the problem (\ref{bs04})}

This section consists of two subsections.  In Subsection 3.1, we consider 
the existence and uniqueness of solutions of (\ref{bs04}).  In 
Subsection 3.2, we obtain some regularity of solutions of (\ref{bs04}).  
To establish the results in these subsections, we use the penalized 
problem for (\ref{bs04}):  
\begin{equation}
\label{penalty}
   \left\{\begin{array}{ll}
       \displaystyle{u^\e_t-u^\e_{xx}+\be u^\e
      +\zeta_\e(u^\e-\varphi)=0} & 
         {\rm in\ }(0,T)\times\om, \\
       u^\e(0,x)=\varphi(x) & {\rm for\ }x\in\oom, \\
       u^\e(t,\pm 1)=\varphi(\pm 1) 
       & {\rm for\ }t\in (0,T),  
          \end{array}\right.
\end{equation}
where $ \e>0 $, $ \zeta_\e(r):=\zeta(r/\e) $ and $ \zeta $ is a smooth 
function such that 
$$
  \zeta^\prime\geq 0,\ \zeta^{\prime\prime}\leq 0\ 
{\rm on\ }\bR,\ 
$$
\subsection{Existence and uniqueness of solutions}

In this subsection we prove the following theorem.  
\begin{theorem}
\label{th301}
There exists a unique solution $ u $ of {\rm (}\ref{bs04}{\rm )} 
in the a.e.\;sense and in the sense of viscosity solutions such that 
$ u\in W^{1,2,2}((0,T)\times\om)\cap C([0,T)\times\oom) $.  
\end{theorem}
\begin{remark}
\label{rem301}{\rm 
See \cite{cr;is;li;92} or \cite{ko;04} for the definition and the theory 
of viscosity solutions.  
}
\end{remark}

By \cite[Chapter 2,\ 3]{ben;li;78}, there is a unique solution 
$ u^\e $ of (\ref{penalty}) in the sense that 
$ u^\e-\varphi\in L^2(0,T;H_0^1(\om)) $, 
$ u^\e_t\in L^2(0,T;L^2(\om)) $ and 
$$
   \int_\om \{u^\e_t\phi+u^\e_x\phi_x+\be u^\e\phi
       +\zeta_\e(u^\e-\varphi)\phi\}dx=0 \quad
    {\rm for\ all}\ \phi\in H_0^1(\om)\ {\rm and\ a.e.\;}t\in(0,T).  
$$
In addition, $ u^\e\in L^2(0,T;H^2(\om))\cap C([0,T);C(\oom)) $.  

We derive some estimates of $ u^\e $.  By the maximum principle, we get 
\begin{equation}
\label{esti307}
    \sup_{\e>0}\|u^\e\|_{C([0,T)\times\oom)}
   \leq \|\varphi\|_{C(\oom)}.  
\end{equation}
By the same arguments as in \cite[Chapter 2, Section 2.4]{ben;li;78} we 
have 
\begin{equation}
\label{esti301}
   \sup_{\e>0}(\|u^\e\|_{L^\infty(0,T;H^1(\om))}
     +\|u^\e_t\|_{L^2(0,T;L^2(\om))})<+\infty.
\end{equation}
To estimate $ \zeta_\e(u^\e-\varphi) $ and $ u^\e_{xx} $, we need 
the following lemma.  
\begin{lemma}
\label{lem309}
There is $ M_1>0 $ such that $ u^\e\geq \varphi-M_1\e $ in 
$ [0,T)\times\oom $ for all $ \e>0 $.  
\end{lemma}

\noindent
{\bf Proof.}  Since $ \varphi $ is Lipschitz continuous and convex in 
$ \om $, we can show that 
\begin{equation}
\label{ineq310}
   \int_\om\varphi_x(x)\phi_x(x)dx\leq 0\quad{\rm for\ all\ }
   \phi\in H_0^1(\om)\ {\rm satisfying\ }\phi\geq 0\ {\rm in\ }\om.  
\end{equation}

Set $ \uu(t,x):=\varphi(x)-M_1\e $.  Then, we use the above inequality 
to obtain 
$$
   \int_\om \{\uu_t\phi+\uu_x\phi_x+\be \uu\phi+\zeta_\e(\uu-\varphi)\phi\}dx
    \leq \int_\om\{\be\varphi+\zeta(-M_1)\}\phi dx
$$
for any $ t\in (0,T) $ and $ \phi\in H_0^1(\om) $ satisfying 
$ \phi\geq 0 $ in $ \om $.  Taking $ M_1>0 $ such that 
$ \be\|\varphi\|_{C(\oom)}+\zeta(-M_1)\leq 0 $, we easily see that 
$ \uu $ is a weak subsolution of (\ref{penalty}).  Hence we have 
$ \uu\leq u^\e $ in $ [0,T)\times\om $ 
by the maximum principle.  Therefore we obtain the result.  $ \square $

\bigskip
\noindent
Hence from (\ref{esti301}) and Lemma \ref{lem309}, we get 
\begin{eqnarray}
\label{esti305}
   & & \sup_{\e>0}
  \|\zeta_\e(u^\e-\varphi)\|_{C([0,T)\times\oom)}<+\infty,\ 
  \sup_{\e>0}\|u^\e_{xx}\|_{L^2((0,T)\times\om)}<+\infty.
\end{eqnarray}

\bigskip
\noindent
{\bf Proof of Theorem \ref{th301}.}  By (\ref{esti307}), 
(\ref{esti301}), (\ref{esti305}) and Sobolev imbedding, we can extract 
a subsequence $ \{\e_n\}_{n=1}^{+\infty} $, $ \e_n\searrow 0 $ such 
that for any $ T^\prime\in (0,T) $ and $ \lambda\in (0,1/2) $, as 
$ n\to +\infty $, 
\begin{eqnarray}
\label{conv301}
  & &  u^{\e_n}\longrightarrow u 
    \quad{\rm in\ } 
    C^{\lambda/2,\lambda}([0,T^\prime]\times\oom),\\ 
\label{conv302}
  & & (u^{\e_n}_t,u^{\e_n}_x,u^{\e_n}_{xx}) \longrightarrow 
    (u_t,u_x,u_{xx})\quad{\rm weakly\ in\ } (L^{2}((0,T)\times\om))^3.  
\end{eqnarray}
We can see that $ u $ is a unique solution of (\ref{bs04}) in the 
a.e.\;sense and in the viscosity sense (cf. \cite[Chapter 3]{ben;li;78}, 
\cite{cr;is;li;92} and \cite{ko;04}).  Thus we complete the proof.  
$ \square $

\subsection{Some properties of solutions}

The main results of this subsection are stated as follows.  Let $ u $ 
be the solution of {\rm(}\ref{bs04}{\rm)}.  
\begin{theorem}
\label{th303}
Assume $ q<r<qe $ and let $ x_0 $ be given in {\rm (}\ref{free}{\rm )}.  
Then $ u\in W_{loc}^{1,2,\infty}((0,T)\times\om) $ and it satisfies 
the following estimates.  
\begin{itemize}
\item[{\rm (1)}] For any small $ x_1>0 $, there is  $ L_1>0 $ 
such that 
\begin{eqnarray}
\label{esti314}
    & & \|u_t(t,\cdot)\|_{L^\infty(-x_1,x_1)}\leq \frac{L_1}{\sqrt{t}},\ 
      \|u_t(t,\cdot)\|_{L^\infty(\om\backslash(-x_1,x_1))}\leq L_1
    \quad for\ a.e.\; t\in (0,T), \\ 
\label{esti308}
    & & \|u_{xx}(t,\cdot)\|_{L^\infty(-x_1,x_1)}\leq \frac{L_1}{\sqrt{t}},\ 
      \|u_{xx}(t,\cdot)\|_{L^\infty(\om\backslash(-x_1,x_1))}\leq L_1
    \quad for\ a.e.\; t\in (0,T).  
\end{eqnarray}
\item[{\rm (2)}] There is $ L_2>0 $ such that 
$$
|u(t,x)-u(s,y)|\leq L_2(|t-s|^{1/2}+|x-y|) 
     \quad {for\ all\ }(t,x),\ (s,y)\in [0,T)\times\oom.  
$$
\end{itemize}
\end{theorem}
\begin{theorem}
\label{th305}
The $ u_t $ is nonnegative and continuous in $ (0,T)\times\oom $.  
\end{theorem}
\begin{remark}
\label{rem302}{\rm  (1) 
In \cite{jai;lam;lap;90}, similar results to Theorem \ref{th303} are 
obtained in the case $ \om=\bR^N $.  

(2) Theorem \ref{th305} is similar to \cite[Lemma 5]{vm;74} and 
\cite[Corollary 4.2]{fr;75}.  As seen in Section 5, it plays an 
important role to prove Theorem \ref{th201}.  
}
\end{remark}

We prepare some pointwise estimates of solutions of (\ref{penalty}) 
to prove Theorem \ref{th303}.  Let $ u^\e $ be the solution of 
{\rm (}\ref{penalty}{\rm )}.  
\begin{proposition}
\label{prop307}
We obtain 
$$
   \sup_{\e>0}\|u^\e_x\|_{L^\infty((0,T)\times\om)}<+\infty.  
$$
\end{proposition}

\noindent
{\bf Proof.}  
The barrier construction argument yields that 
$ |u^\e_x(t,\pm 1)|\leq\|\varphi\|_{L^\infty(\om)} $ for all 
$ t\in [0,T) $ and $ \e>0 $.  We obtain the result by combining the 
comparison argument for viscosity solutions (cf. 
\cite[Section 7]{is;li;90}) with Lipschitz continuity of $ \varphi $ 
and this estimate.  $ \square $

\begin{lemma}
\label{lem301}
Assume $ q<r<qe $ and let $ x_0 $ be defined by {\rm (}\ref{free}{\rm )}.  
Then there exist $ x_2\in (0,x_0) $ such that 
$ u^\e(t,x)>\varphi(x) $ for all $ t\in (0,T) $, 
$ x\in (-1,x_2) $ and $ \e>0 $.  
\end{lemma}
We can formally show this lemma, according to \cite[p.124]{wi;ho;de;95}.  
Let $ x_0 $ be defined in (\ref{free}) and assume $ x_0\in (0,1) $.  
It is seen that for small $ \e>0 $, 
$$
   u^\e_t(t,x)\longrightarrow u^\e_t(0,x)=\varphi_{xx}(x)-\be\varphi(x)
   \quad{\rm for\ all\ }x\in (0,1)\ {\rm as\ }t\searrow 0.
$$
It follows from the definition of $ \varphi $ that 
\begin{equation}
\label{ineq304}
   \varphi_{xx}(x)-\be\varphi(x)
   =e^{\al x}(-qe^x+r)
     \left\{\begin{array}{ll}
             >0 & {\rm if\ }x\in (0,x_0), \\
             =0  & {\rm if\ }x=x_0,\\
             <0  & {\rm if\ }x\in (x_0,1], 
     \end{array}\right.
\end{equation}
Thus Lemma \ref{lem301} formally holds with $ x_1=x_0 $.  

\bigskip
\noindent
{\bf Proof of Lemma \ref{lem301}.} {\it Step 1.} Let $ \uu $ be the 
classical solution of 
\begin{equation}
\label{para301}
   \left\{\begin{array}{ll}
       \uu_t-\uu_{xx}+\be \uu=0 & 
         {\rm in\ }(0,+\infty)\times\om, \\
       \uu(0,x)=\varphi(x) & {\rm for\ }x\in\oom, \\
       \uu(t,\pm 1)=\varphi(\pm 1) & {\rm for\ }t\in (0,+\infty).  
          \end{array}\right.
\end{equation}
We prove that 
\begin{equation}
\label{esti310}
   u^\e(t,x)\geq \uu(t,x)>\varphi(x)\quad{\rm for\ all\ }t>0,\ 
   x\in (-1,0]\ {\rm and\ }\e>0.  
\end{equation}

Since $ \uu $ is a subsolution of (\ref{penalty}), it follows from 
the maximum principle that $ \uu\leq u^\e $ in $ [0,T)\times\oom $.  
Hence we get (\ref{esti310}) by $ \uu>0 $ in $ (0,+\infty)\times \om $ 
and $ \varphi\equiv 0 $ on $ [-1,0] $.  

{\it Step 2.}  We show that there exist $ t_1 $, $ M_2>0 $ such 
that 
\begin{equation}
\label{esti302}
  \uu(t,x)\geq \varphi(x)+M_2t
     +\int_0^t\frac{e^{-\be s-x^2/4s}}{\sqrt{8\pi s}}ds 
   \quad{\rm for\ all\ }
     t\in (0,t_1]\ {\rm and\ }x\in [0,2x_0/3].  
\end{equation}

For $ a>0 $, define 
\begin{eqnarray}
\label{fund}
  & & E^a(t,x,y):=\frac{1}{\sqrt{4\pi t}}
   \sum_{n=-\infty}^{+\infty}\{e^{-(x-y+4an)^2/4t}
   -e^{-(x+y+4an+2a)^2/4t}\}, \\
\label{eq301}
   & & E_0(t,x,y):=\frac{1}{\sqrt{4\pi t}}e^{-(x-y)^2/4t},\ 
   E^a_1(t,x,y):=E^a(t,x,y)-E_0(t,x,y). 
\end{eqnarray}
Put $ E=E^1 $ and $ E_1=E^1_1 $ for simplicity. 
Then $ \uu $ is given by 
$$
\uu(t,x)=e^{-\be t}\int_0^1 E(t,x,y)\varphi(y)dy
      -\varphi(1)\int_0^t e^{-\be (t-s)}E_y(t-s,x,1)ds. 
$$
Differentiating this formula with respect to $ t $, we have 
$$
   \uu_t(t,x)
   =e^{-\be t}\int_0^1 E_t(t,x,y)\varphi(y)dy
      -\be e^{-\be t}\int_0^1 E(t,x,y)\varphi(y)dy 
        -e^{-\be t}\varphi(1)E_y(t,x,1).  
$$
We use the facts $ E_t=E_{xx}=E_{yy} $, $ \varphi(0)=0 $ and 
the integration by parts to obtain 
\begin{eqnarray*}
   \uu_t(t,x)&=&e^{-\be t}\int_0^1 E(t,x,y)(\varphi_{yy}(y)-\be\varphi(y))dy
      +e^{-\be t}E(t,x,+0)\varphi_y(+0) \\
     &=:&I_{1,1}+I_{1,2}.  
\end{eqnarray*}

We estimate the right-hand side (RHS for short) of the above formula 
to have (\ref{esti302}).  Some calculations yield that for small 
$ t>0 $, $ x\in [-2x_0/3,2x_0/3] $ and 
$ y\in (0,1) $,   
$$
   |E(t,x,y)-E_0(t,x,y)|\leq Ce^{-x_0^2/64t}.  
$$
We observe from $ q<r<qe $ and this estimate that 
\begin{eqnarray*}
   I_{1,1}&\geq &\left(\int_0^{3x_0/4}+\int_{3x_0/4}^1\right)
      E_0(t,x,y)e^{-\be t+\al y}(-qe^y+r)dy-Ce^{-x_0^2/64t} \\
        &\geq & \frac{e^{-\be t-|\al|}}{\sqrt{4\pi t}}
          \int_0^{3x_0/4}e^{-(y-x)^2/4t}(-qe^y+r)dy 
             -Ce^{-x_0^2/64t}\left(\frac{1}{\sqrt{4\pi t}}
             +1\right)  
\end{eqnarray*}
for small $ t>0 $ and $ x\in [-2x_0/3,2x_0/3] $.  By (\ref{ineq304}) 
we have $ -qe^y+r\geq M_{2,1} $ for all $ y\in [0,3x_0/4] $ and some 
$ M_{2,1}>0 $.  Thus 
$$
   I_{1,1}\geq M_2
   \quad{\rm for\ small\ }t>0,\ x\in [-2x_0/3,2x_0/3]\ 
   {\rm and\ some\ }M_2>0.  
$$
Since it is seen that $ I_{1,2}\geq e^{-\be t-x^2/4t}/\sqrt{8\pi t} $ 
for small $ t>0 $ and all $ x\in\oom $, we obtain 
$$
    u_t(t,x)\geq M_2+\frac{e^{-\be t-x^2/4t}}{\sqrt{8\pi t}}
   \quad{\rm for\ small\ }t\ {\rm and\ }x\in [0,2x_0/3].  
$$
Therefore, for sufficiently small $ t_1>0 $, we have (\ref{esti302}) 
by integrating both sides of this inequality on $ [0,t] $ for 
all $ t\in (0,t_1) $.  

{\it Step 3.}  Set 
$ 
    \ox(t):=\sup\{y\in[0,1]\;|\;\uu(t,x)>\varphi(x)\ {\rm for\ all\ }
      x\in [0,y)\}  
$
for each $ t\in[t_1,T] $ and define 
$ \ds{x_{2,1}:=\inf_{t\in[t_1,T]}\ox(t)} $. We claim $ x_{2,1}>0 $.  

Suppose $ x_{2,1}=0 $.  Then for each 
$ n\in\bN $, there exists 
$ t_n\in [t_1,T] $ such that $ \ox(t_n)\leq 1/n $.  Extracting 
a subsequence if necessary, we may assume 
$ t_n\longrightarrow \widetilde t\in [t_1,T] $ as $ n\to+\infty $.  
Noting that $ \uu(t,\ox(t))=\varphi(\ox(t)) $, we easily see that 
$$
    \uu(\widetilde t,0)=\lim_{n\to+\infty}\uu(t_n,\ox(t_n))
      =\lim_{n\to+\infty}\varphi(\ox(t_n))=\varphi(0)=0.  
$$
This contradicts to (\ref{esti310}).  Hence the claim of this step is 
proved.  

Putting $ x_2:=\min\{2x_0/3,x_{2,1}\} $,  we obtain the desired 
result.  $ \square $

\begin{remark}
\label{rem303}{\rm 
It readily follows from Lemma \ref{lem301} that 
$ \zeta_\e(u^\e-\varphi)\equiv 0 $ in $ [0,T)\times(-x_2,x_2) $ for all 
$ \e>0 $.  Hence the boot-strap argument yields that 
$ u^\e\in C^\infty((0,T)\times\om) $. 
}
\end{remark}

Based on (\ref{esti305}) and Lemma \ref{lem301}, we prove the 
following theorem.  
\begin{theorem}
\label{th308}
Assume $ q<r<qe $.  Let $ x_0 $ be given in {\rm (}\ref{free}{\rm )} 
and $ x_2\in (0,x_0) $ in Lemma \ref{lem301}.  Then, for each 
$ x_3\in (0,x_2) $, there exists $ L_3>0 $ such that for any $ \e>0 $, 
\begin{eqnarray}
\label{esti316}
  & &  \|u^\e_t(t,\cdot)\|_{L^\infty(-x_3,x_3)}\leq \frac{L_3}{\sqrt{t}},\ 
      \|u^\e_t(t,\cdot)\|_{L^\infty(\om\backslash(-x_3,x_3))}\leq L_3\quad
    {\ for\ all\ }t\in (0,T),\\
\label{esti317}
  & &  \|u^\e_{xx}(t,\cdot)\|_{L^\infty(-x_3,x_3)}\leq \frac{L_3}{\sqrt{t}},\ 
    \|u^\e_{xx}(t,\cdot)\|_{L^\infty(\om\backslash(-x_3,x_3))}\leq L_3\quad
    {\ for\ all\ }t\in (0,T).  
\end{eqnarray}
\end{theorem}

\noindent
{\bf Proof.}  The $ u^\e $ is given by 
$$
  u^\e(t,x)=\uu(t,x)
    -\int_0^t \kern -1.5mm\int_\om e^{-\be(t-s)}E(t-s,x,y)
    \zeta_\e(u^\e(s,y)-\varphi(y))dyds,  
$$
where $ \uu $ is the solution of (\ref{para301}) and $ E $ is defined 
by (\ref{fund}).  Fix $ x_3\in (0,x_2/2) $.  We divide our 
consideration into three cases.  

{\it Case 1.} $ |x|<x_3 $.  

Differentiating $ u^\e $ with respect to $ t $, we get 
\begin{eqnarray}
\label{eq303}
  u_t^\e(t,x)&=& \uu_t(t,x)-\zeta_\e(u^\e(t,x)-\varphi(x)) \\
\nonumber
   & & -\int_0^t \kern -1.5mm\int_\om e^{-\be(t-s)}E_t(t-s,x,y)
    \zeta_\e(u^\e(s,y)-\varphi(y))dyds \\
\nonumber
   & & +\be\int_0^t \kern -1.5mm\int_\om e^{-\be(t-s)}E(t-s,x,y)
    \zeta_\e(u^\e(s,y)-\varphi(y))dyds \\
\nonumber
   &=:& I_{2,1}+I_{2,2}+I_{2,3}+I_{2,4}.  
\end{eqnarray}
It follows from the standard theory for parabolic equations that 
$ |I_{2,1}|\leq C/\sqrt{t} $ for all $ t\in (0,T) $ and $ x\in \oom $.  
Besides, from (\ref{esti305}) we get $ |I_{2,2}+I_{2,4}|\leq C $ for 
all $ t\in (0,T) $ and $ x\in\oom $.  As for $ I_{2,3} $, noting that 
$ (x-y)^2\geq x_3^2 $ for all $ x\in [-x_3,x_3] $ and 
$ y\in [2x_3,1) $, we observe by (\ref{esti305}) and Lemma \ref{lem301} 
that 
\begin{eqnarray*}
   |I_{2,3}|
&\leq & C\left(\left|\int_0^t \kern -1.5mm\int_{x_2}^1E_{0,t}(t-s,x,y)
     dy\right|
       +\left|\int_0^t \kern -1.5mm\int_{x_2}^1
    E_{1,t}(t-s,x,y)dy\right|\right)\\
   & \leq & C\left\{\int_0^t(t-s)^{-3/2}e^{-x_3^2/4(t-s)}ds+1\right\}
      \leq M_{3,1}.  
\end{eqnarray*}
for all $ t\in (0,T) $, $ x\in (-x_3,x_3) $ and $ \e>0 $.  Here and in 
the sequel, $ M_{3,i} $'s $ (i\geq 1) $ are constants depending on 
$ x_3 $.  Consequently we obtain  
$$
   \|u_t^\e(t,\cdot)\|_{L^\infty(-x_3,x_3)}\leq\frac{M_{3,2}}{\sqrt{t}}\quad 
     {\rm for\ all\ }t\in (0,T)\ {\rm and\ small\ }\e>0.  
$$

{\it Case 2.} $ x_3\leq |x|\leq 1 $.  

Assume that $ x_3\leq x\leq 1 $.  By using (\ref{eq303}) and 
(\ref{esti305}), it is seen that 
\begin{eqnarray*}
   |u_t^\e(t,x_3)|&\leq & |\uu_t(t,x_3)|
     +C\left|\int_0^t \kern -1.5mm\int_{x_3}^1E_t(t-s,x_3,y)
    dyds\right| +C \\
   &=:&I_{3,1}+I_{3,2}+C\quad{\rm for\ all\ }t\in (0,T)
   \ {\rm and\ small\ }\e>0.  
\end{eqnarray*}

We estimate $ I_{3,1} $ and $ I_{3,2} $.  We observe by the integration 
by parts that 
\begin{eqnarray*}
  |I_{3,1}|&\leq & \left|\int_0^1E_y(t,x_3,y)\varphi_y(y)dy\right|
            +\be\|\varphi\|_{L^\infty(\om)} \\
     &=&\left|\int_0^{x_3/2}E_y(t,x_3,y)\varphi_y(y)dy\right|
       +\left|\int_{x_3/2}^{1}E_y(t,x_3,y)\varphi_y(y)dy\right|+C\\
     &=:& I_{3,1,1}+I_{3,1,2}.  
\end{eqnarray*}
Using $ (x_3-y)^2\geq x_3^2/4 $ for all $ y\in (0,x_3/2] $, 
we can estimate 
$ |I_{3,1,1}|\leq Ct^{-3/2}e^{-x_3^2/16t}\leq M_{3,3} $ for all 
$ t\in (0,T) $ and $ \e>0 $.  Since $ \varphi $ is smooth in $ (0,1) $, 
it follows from the integration by parts and the fact 
$ |E(t,x_3,x_3/2)|+|E(t,x_3,1)|\leq Ct^{-1/2}e^{-x_3^2/16t} $ 
that for all $ t\in (0,T) $ and $ \e>0 $,   
$$
  I_{3,1,2}\leq |E(t,x_3,1)\varphi_y(1)-E(t,x_3,x_3/2)\varphi_y(x_3/2)| 
     +\left|\int_{x_3/2}^{1}E(t,x_3,y)\varphi_{yy}(y)dy\right|
    \leq  M_{3,4}.  
$$
Hence $ |I_{3,1}|\leq M_{3,5} $ for some $ M_{3,5}>0 $.  On the other 
hand, it is observed by the same argument as the estimate for 
$ I_{2,3} $ in Case 1 that 
$$
   |I_{3,2}|\leq C\int_0^t(t-s)^{-3/2}e^{-x_3^2/16(t-s)}ds
    \leq M_{3,6}\quad {\rm for\ all\ }t\in (0,T)\ 
    {\rm and\ }\e>0.  
$$
Therefore we conclude that $ |u_t^\e(t,x_3)|\leq M_{3,7} $ for all 
$ t\in (0,T) $ and $ \e>0 $.  

We provide an estimate for $ \|u^\e_t(t,\cdot)\|_{L^\infty(x_3,1)} $.  
Differentiating the equation of (\ref{penalty}) with respect to 
$ t $, we have 
$$
   \left\{\begin{array}{ll}
      U^\e_t-U^\e_{xx}+\be U^\e+\zeta_\e^\prime(u^\e-\varphi)U^\e=0 & 
         {\rm in\ }(0,T)\times(x_3,1), \\
     U^\e(0,x)=\varphi_{xx}(x)-\be\varphi(x) & 
       {\rm for\ }x\in(x_3,1), \\
     |U^\e(t,x_3)|\leq M_{3,7},\ U^\e(t,1)=0 
   & {\rm for\ }t\in(0,T), 
          \end{array}\right.
$$
where $ U^\e:=u_t^\e $.  Noting that $ |U^\e(0,\cdot)|\leq C $ on 
$ [x_3,1] $ and that $ \zeta_\e^\prime\geq 0 $, we obtain from the 
maximum principle $ \|u_t^\e\|_{L^\infty(x_3,1)}\leq M_{3,8} $ for all 
$ t\in(0,T) $, small $ \e>0 $.  

We can get $ \|u_t^\e(t,\cdot)\|_{L^\infty(-1,-x_3)}\leq M_{3,9} $ for 
all $ t\in(0,T) $ and small $ \e>0 $ by the same way as above.  

Consequently, for each $ x_3\in (0,x_2/2) $, there is 
$ L_3\geq\max\{M_{3,2},M_{3,8},M_{3,9}\} $ such that (\ref{esti316}) 
holds for all $ t\in(0,T) $.  The (\ref{esti317}) follows from 
(\ref{esti307}), (\ref{esti305}) and (\ref{esti316}).  
$ \square $

\bigskip
\noindent
{\bf Proof of Theorem \ref{th303}.} Proposition \ref{prop307} and 
Theorem \ref{th308} yield that there is a subsequence 
$ \{\e_n\}_{n=1}^{+\infty} $, $ \e_n\searrow 0 $, such that as 
$ n\to+\infty $, 
$$
   (u_x^{\e_n},\sqrt{t}u_t^{\e_n},
 \sqrt{t}u_{xx}^{\e_n})\longrightarrow (U_1,U_2,U_3)\quad
   {\rm weakly\ star\ in\ }(L^\infty((0,T)\times\om))^3.
$$
By $ L^2((0,T)\times\om)\subset L^1((0,T)\times\om) $, we can use 
(\ref{conv302}) to have 
$ (U_1,U_2,U_3)=(u_x,\sqrt{t}u_t,\sqrt{t}u_{xx}) $.  
Hence the $ u\in W^{1,2,\infty}_{loc}((0,T)\times\om) $ follows from 
(\ref{esti307}) and these convergences.  

Set $ x_1\in(0,x_3) $.  The (\ref{esti314}) and (\ref{esti308}) are 
derived from the above convergences and Theorem \ref{th308}.  The 
asserion of (2) follows from (\ref{conv301}), Theorem \ref{th308} and 
Proposition \ref{prop307}.  $ \square $ 

\bigskip
We prepare some estimates for $ u_{xt} $ and $ u_{tt} $ to show Theorem 
\ref{th305}.  
\begin{proposition}
\label{prop303}
There exists $ C>0 $ such that for any small $ \sigma>0 $, 
$$
   \sup_{t\in(\sigma,T)}\|u_{xt}(t,\cdot)\|_{L^2(\om)}
   +\|u_{tt}\|_{L^2((\sigma,T)\times\om)}\leq \frac{C}{\sigma}.  
$$
\end{proposition}

\noindent
{\bf Proof.}  Let $ u^\e $ be the solution of (\ref{penalty}).  Then 
$ U^\e:=u^\e_t $ satisfies 
$$
   \left\{\begin{array}{ll}
     U^\e_t-U^\e_{xx}+\be U^\e+\zeta_\e^\prime(u^\e-\varphi)U^\e=0 
        & {\rm in\ }(\sigma,T)\times\om, \\
     \ds{|U^\e(\sigma,x)|\leq \frac{C}{\sqrt{\sigma}}} & {\rm for\ }x\in\oom,\\
     U^\e(t,\pm 1)=0 & {\rm for\ }t\in (\sigma,T).  
   \end{array}\right.  
$$
Here the second inequality follows from (\ref{esti314}).  The same 
argument as in the proof of \cite[Lemma 3.4]{fr;ki;75} yields that 
$$
   \sup_{t\in(\sigma,T)}\|u^\e_{xt}(t,\cdot)\|_{L^2(\om)}
   +\|u^\e_{tt}\|_{L^2((\sigma,T)\times\om)}\leq \frac{C}{\sigma}.  
$$
Sending $ \e\to 0 $, we have the result.  $ \square $

\bigskip
\noindent
{\bf Proof of Theorem \ref{th305}.}  
{\it Step 1.} We claim that $ u_t $ is continuous in $ (0,T)\times\om $.  

We observe from the regularity theory for parabolic equations that 
$ u_t $ is continuous in $ \{(t,x)\;|\;0<t<T,\ u(t,x)>\varphi(x)\} $.  
It is obvious that $ u_t $ is so in $ \{(t,x)\;|\;0<t<T,\ x^*(t)<x<1\} $.  
The continuity of $ u_t $ in $ \{(t,x^*(t))\;|\;0<t<T\} $ can be 
proved by Proposition \ref{prop303} and the same argument as the proof 
of \cite[Corollary 4.2]{fr;75}.  Hence we have the claim.  

{\it Step 2.} We show $ u_t\geq 0 $ in $ (0,T)\times\oom $.     

We modify (\ref{penalty}) as follows.  Let $ \{\varphi_\delta\}_{\delta>0} $ 
be a sequence of $ C^2 $ and convex functions satisfying 
$ \|\varphi_\delta-\varphi\|_{W^{1,\infty}(\om)}\longrightarrow 0 $ 
as $ \delta\to 0 $.  We consider the following instead of (\ref{penalty}).  
\begin{equation}
\label{penalty03}
   \left\{\begin{array}{ll}
      u^{\delta,\e}_t-u^{\delta,\e}_{xx}+\be u^{\delta,\e}
   +\zeta_\e(u^{\delta,\e}-\varphi_\delta-M_4\e)=0 & 
         {\rm in\ }(0,T)\times\om, \\
     u^{\delta,\e}(0,x)=\varphi_\delta(x) & {\rm for\ }x\in\oom, \\
     u^{\delta,\e}(t,\pm 1)=\varphi_\delta(\pm 1) & {\rm for\ }t\in(0,T)
          \end{array}\right.
\end{equation}
Here $ M_4>0 $ is chosen so that 
$ \ds{\be\sup_{\delta\in(0,1)}\|\whvphi_{\delta}\|_{L^\infty(\om)}
+\zeta(-M_4)\leq 0} $.  Then, there is a unique classical solution 
$ u^{\delta,\e} $ of (\ref{penalty03}) and it satisfies 
\begin{equation}
\label{conv306}
   \lim_{\delta\to 0}\lim_{\e\to 0}u^{\delta,\e}=u\quad
   {\rm locally\ uniformly\ in\ }[0,T)\times\oom. 
\end{equation}

We differentiate (\ref{penalty03}) with respect to $ t $ and set 
$ U^{\delta,\e}:=u_t^{\delta,\e} $.  Then we have 
$$
   \left\{\begin{array}{ll}
      U^{\delta,\e}_t-U^{\delta,\e}_{xx}+\be U^{\delta,\e}
   +\zeta_\e^\prime(u^{\delta,\e}-\varphi_\delta-M_4\e)U^{\delta,\e}=0 & 
         {\rm in\ }(0,T)\times\om, \\
     U^{\delta,\e}(0,x)=\varphi_{\delta,xx}-\be \whvphi_{\delta}
      -\zeta(-M_4) & {\rm for\ }x\in\oom, \\
     U^{\delta,\e}(t,\pm 1)=0 & {\rm for\ }t\in(0,T).
          \end{array}\right.
$$
The $ U^{\delta,\e}(0,\cdot)\geq 0 $ on $ \oom $ follows from 
$ \varphi_{\delta,xx}(x)\geq 0 $ and the choice 
of $ M_4 $.  
Hence we apply the maximum principle to obtain $ U^{\delta,\e}\geq 0 $ 
in $ [0,T)\times\oom $.  The (\ref{conv306}) and this result yield 
that $ u(\cdot,x) $ is nondecreasing for each $ x\in\oom $.  Hence 
we have the result.  $ \square $

\section{The discrete Morse semiflow}

Fix $ h>0 $.  As briefly mentioned in Section 2, for each 
$ m=1,2,\ldots, [T/h] $, there is a unique minimizer $ u_m\in\cK $ 
of the functional $ J_m $.  Moreover, $ u_m $ satisfies (\ref{ellvi}) 
in the weak sense (cf. \cite[Chapter 3]{ben;li;78}).  We call the 
sequence $ \{u_m\}_{m=0}^{[T/h]} $ the DMS-BS.  This section is 
devoted to some properties of the DMS-BS and to its free boundary.  

In Subsection 4.1, we discuss some properties 
of the DMS-BS.  To prove Theorem \ref{th402} of this subsection, 
we need an estimate of the difference $ u_m-u_{m-1} $, 
Theorem \ref{th406}.  Since its proof consists of lengthy and careful 
calculations, it is given 
in Subsection 4.3.  In Subsection 4.2, we consider the existence and 
uniquess of the free boundary of the DMS-BS.  

\subsection{Some properties of the DMS-BS}

First, we have the monotone property of the DMS-BS.  
\begin{theorem}
\label{mono}
Let $ \{u_m\}_{m=0}^{[T/h]} $ be the DMS-BS.  Then 
$ u_{m-1}\leq u_m $ on $ \oom $ for all $ m=1,2,\ldots,[T/h] $ and 
$ h>0 $.  
\end{theorem}
This theorem can be easily proved by the maximum principle 
and induction.  Hence we omit the proof.  

In the following part of this subsection, we show the time-discrete 
analogues to Theorems \ref{th301} and \ref{th303}.  
\begin{theorem}
\label{th401}
For each $ h>0 $ and $ m=1,2,\ldots,[T/h] $, $ u_m $ is a unique solution 
of {\rm (}\ref{ellvi}{\rm )} in the a.e.\;sense and in the viscosity sense.  
In addition, $ \{u_m\}_{m=0}^{[T/h]} $ satisfies 
$$
   \sup_{h>0}\left(h\sum_{m=1}^{[T/h]}\left\|\frac{u_m-u_{m-1}}{h}
   \right\|_{L^2(\om)}^2+\sup_{0\leq m\leq[T/h]}\|u_m\|_{H^1(\om)}
    +h\sum_{m=1}^{[T/h]}\|u_{m,xx}\|_{L^2(\om)}^2
   \right)<+\infty.  
$$
\end{theorem}
\begin{theorem}
\label{th402}
Assume $ q<r<qe $ and let $ x_0 $ be defined in {\rm(}\ref{free}{\rm)}.  
Then there exist $ x_4\in (0,x_0) $ and $ h_1>0 $ such that for each 
$ h\in(0,h_1) $, 
$ \{u_m\}_{m=0}^{[T/h]}\subset W^{2,\infty}(\om) $ and it satisfies 
the following estimates.  
\begin{itemize}
\item[\rm (1)] There are $ L_4 $, $ L_5>0 $ depending on 
$ x_4 $ such that for each $ h\in (0,h_1) $ and 
$ m=1,2,\ldots,[T/h] $, 
\begin{eqnarray}
\label{esti401}
  & & \frac{|u_m(x)-u_{m-1}(x)|}{h}\leq \left\{\begin{array}{ll}
         \ds{L_4\max\left\{\frac{1}{\sqrt{mh}},\sqrt{\frac{|\log h|}{L_5}}
               \right\}}
          & {for\ }|x|<x_4,\\
         \ds{L_4\sqrt{\frac{|\log h|}{L_5}}}
          & {for\ }|x|\geq x_4,
      \end{array}\right. \\
\label{esti402}
  & & |u_{m,xx}(x)|\leq\left\{\begin{array}{ll}
         \ds{L_4\max\left\{\frac{1}{\sqrt{mh}},\sqrt{\frac{|\log h|}{L_5}}
               \right\}}
          & {for\ }a.e.\;|x|<x_4,\\
         \ds{L_4\sqrt{\frac{|\log h|}{L_5}}}
          & {for\ }a.e.\;|x|\geq x_4.
      \end{array}\right. 
\end{eqnarray}
\item[\rm (2)] There are $ L_6 $, $ L_7>0 $ depending on $ x_4 $ 
such that 
$$
   |u_m(x)-u_n(y)|\leq \left\{\begin{array}{ll}
      L_6(|(m-n)h|^{1/2}+|x-y|) & {if\ }m,n\leq L_7(h|\log h|)^{-1},\\
        L_6(|(m-n)h|\sqrt{|\log h|}+|x-y|) & {if\ }m,n\geq L_7(h|\log h|)^{-1},
    \end{array}\right.
$$
for all $ m $, $ n=0,1,\ldots,[T/h] $, $ x $, $ y\in\oom $ and 
$ h\in(0,h_1) $. 
\end{itemize}
\end{theorem}
\begin{remark}{\rm 
The $ \sqrt{|\log h|} $ appears in Theorem \ref{th402} by some technical 
reasons.  
}
\end{remark}

To prove Theorems \ref{th401} and \ref{th402}, we introduce the penalized 
problem to (\ref{ellvi}):  Put $ u_0^\e:=\varphi $ and consider 
\begin{equation}
\label{ellvi02}
   \frac{u_m^\e-u^\e_{m-1}}{h}-u^\e_{m,xx}+\be u_m^\e
    +\zeta_\e(u_m^\e-\varphi)=0\quad 
    {\rm in\ }\om,\ u_m^\e(\pm 1)=\varphi(\pm 1).  
\end{equation}
Here $ \zeta_\e $ is the same function as in Section 3.  

We observe that for each $ h>0 $ and $ m=1,2,\ldots,[T/h] $, 
there uniquely exists a weak solution $ u^\e_m $ of (\ref{ellvi02}) 
in the sense that $ u^\e_m-\varphi\in H_0^1(\om) $ and 
\begin{equation}
\label{wsol02}
  \int_\om \left\{\frac{u_m^\e-u^\e_{m-1}}{h}\phi
     +u^\e_{m,x}\phi_x+\be u_m^\e\phi
       +\zeta_\e(u_m^\e-\varphi)\phi\right\}dx=0\ \ 
  {\rm for\ all\ }\phi\in H_0^1(\om).
\end{equation}
In addition, the regularity theory for elliptic equations yields that 
$ u^\e_m\in H^2(\om)\cap C^2(\om) $.  Thus $ u^\e_m $ is a classical 
solution of (\ref{ellvi02}).  

We derive some uniform estimates of $ \{u_m^\e\}_{m,\e} $ to prove 
Theorem \ref{th401}.   We get from the maximum principle and induction 
\begin{equation}
\label{esti403}
   \sup_{h>0,\e>0}
    \left(\sup_{0\leq m\leq [T/h]}\|u_m^\e\|_{C(\oom)}\right)
    \leq \|\varphi\|_{C(\oom)}.  
\end{equation}
By a similar argument to the proof of Lemma \ref{lem309} and induction, 
we have  
\begin{lemma}
\label{lem407}
We have $ u_m^\e\geq \varphi-M_1\e $ on $ \oom $ for all $ \e>0 $, 
$ m=1,2,\ldots,[T/h] $ and $ h>0 $.  
Here $ M_1 $ is the same constant as in Lemma \ref{lem309}.  
\end{lemma}
The following estimate is a time-discrete analogue to 
(\ref{esti301}).  
\begin{proposition}
\label{prop401}
We have 
$$
   \sup_{h>0,\e>0}\left(h\sum_{m=1}^{[T/h]}
    \left\|\frac{u_m^\e-u_{m-1}^\e}{h}\right\|^2_{L^2(\om)}
     +\sup_{0\leq m\leq [T/h]}\|u_m^\e\|_{H^1(\om)}\right)
    <+\infty.
$$
\end{proposition}

\noindent
{\bf Proof.} Put $ \phi:=(u^\e_m-u^\e_{m-1}) $ in (\ref{wsol02}).  Using 
$ \omega\widetilde\omega\leq (\omega^2+\widetilde\omega^2)/2 $ 
($ (\omega,\widetilde\omega)=(u^\e_m,u^\e_{m-1}) $, 
($ u^\e_{m,x},u^\e_{m-1,x}) $), we have 
\begin{eqnarray*}
   & & \int_\om \Bigg[h\left|\frac{u^\e_m-u^\e_{m-1}}{h}\right|^2
     +\frac{(u_{m,x}^\e)^2-(u^\e_{m-1,x})^2}{2}
      +\be \frac{(u_m^\e)^2-(u^\e_{m-1})^2}{2} \\
   & & \hspace{70mm}
    +\zeta_\e(u_m^\e-\varphi)(u^\e_m-u^\e_{m-1})
    \Bigg]dx\leq 0.
\end{eqnarray*}

Since it is easily seen from Lemma \ref{lem407} that 
$$
   \zeta_\e(u_m^\e-\varphi)(u^\e_m-u^\e_{m-1})
   \geq -C|u^\e_m-u^\e_{m-1}|
   \geq -\frac{Ch}{2}-\frac{h}{2}
        \left|\frac{u^\e_m-u^\e_{m-1}}{h}\right|^2,  
$$
we have 
\begin{eqnarray*}
   & & \int_\om \Bigg\{\frac{h}{2}\left|\frac{u^\e_m-u^\e_{m-1}}{h}\right|^2
     +\frac{(u_{m,x}^\e)^2-(u^\e_{m-1,x})^2}{2}
      +\be \frac{(u_m^\e)^2-(u^\e_{m-1})^2}{2}\Bigg\}dx\leq Ch.
\end{eqnarray*}
for all $ \e>0 $, $ m=1,2,\ldots,[t/h] $ and $ h>0 $.  
Summing up these inequalities from $ m=1 $ to $ m=[t/h] $, 
we obtain 
$$
   h\sum_{m=1}^{[t/h]}\left\|\frac{u^\e_m-u^\e_{m-1}}{h}\right\|_{L^2(\om)}^2
    +\min\{1,\be\}\|u^\e_{[t/h]}\|_{H^1(\om)}^2 
   \leq \max\{1,\be\}\|\varphi\|_{H^1(\om)}^2+C.  
$$
Since $ t\in (0,T) $ is arbitrary, we have the result.  $ \square $

\bigskip
From Lemma \ref{lem407} and Proposition \ref{prop401}, we have 
\begin{eqnarray}
\label{esti405}
   & & \sup_{h>0,\e>0}\left(\sup_{0\leq m\leq [T/h]}
  \|\zeta_\e(u_m^\e-\varphi)\|_{C(\oom)}\right)<+\infty,\ 
   \sup_{h>0,\e>0}h\sum_{m=1}^{[T/h]}
    \|u_{m,xx}^\e\|_{L^2(\om)}^2<+\infty.
\end{eqnarray}

\noindent
{\bf Proof of Theorem \ref{th401}.}  
From (\ref{esti403}), Proposition \ref{prop401} and (\ref{esti405}), 
we observe that for each $ h>0 $ and $ m=1,2,\ldots,[T/h] $, 
$ \ds{\sup_{\e>0}\|u_m^\e\|_{H^2(\om)}\leq C/h} $.  Hence 
applying Sobolev imbedding, we can extract a subsequence 
$ \{\e_n\}_{n=1}^{+\infty} $, $ \e_n\searrow 0 $ such that 
as $ n\to+\infty $, 
\begin{equation}
\label{conv401}
  u_m^{\e_n}\longrightarrow \wu_m\quad{\rm in\ }C^1(\oom),\ 
  u_{m,xx}^{\e_n}\longrightarrow \wu_{m,xx}\quad
    {\rm weakly\ in\ }L^2(\om),  
\end{equation}
for all $ m=1,2,\ldots,[T/h] $ and $ h>0$.  Thus $ \wu_m $ is 
a solution of (\ref{ellvi}) in the a.e.\;sense and in the viscosity 
sense.  The $ \wu_m=u_m $ follows from the uniqueness of solutions of 
(\ref{ellvi}). 
The estimates in Theorem \ref{th401} follows from 
Proposition \ref{prop401}, (\ref{esti405}) and (\ref{conv401}).  
$ \square $

\bigskip
We provide some pointwise estimates for $ \{u_m^\e\}_{m,\e} $.  
We can show by a similar argument to the proof of Proposition 
\ref{prop307} that 
$$
   \sup_{h>0,\e>0}\left(\sup_{0\leq m\leq [T/h]}
    \|u_{m,x}^\e\|_{L^\infty(\om)}\right)<+\infty.  
$$
The following theorem plays a crucial role to prove Theorem \ref{th402} 
and Theorem \ref{th404} in Subserction 4.2 below.  
\begin{theorem}
\label{th405}
Assume $ q<r<qe $ and let $ x_0 $ be given in {\rm (}\ref{free}{\rm )}.  
Then there exists $ h_2>0 $ such that 
$ u_m^\e>\varphi $ in $ (-1,x_0+\sqrt{h}/2) $ for all $ \e\in (0,h^4) $, 
$ m=1,2,\ldots,[T/h] $ and $ h\in (0,h_2) $.  
\end{theorem}
To prove Theorem \ref{th405}, we prepare some lemmas.  
\begin{lemma}
\label{lem401}
Let $ M_1>0 $ be given in Lemma \ref{lem309}.  
Then for each $ h>0 $, $ m=1,2,\ldots,[T/h] $ and 
$ \e\in (0,h^4) $, we have $ u_m^\e\geq u_{m-1}^\e-M_1h^4 $ on 
$ \oom $.  
\end{lemma}
This lemma is a substitute for Theorem \ref{mono}.  Because we do not 
know such a monotone property for $ \{u_m^\e\}_{m=0}^{[T/h]} $ as Theorem 
\ref{mono} holds since it may happen $ u_m^\e(x)<\varphi(x) $ 
in view of Lemma \ref{lem409}.  

\bigskip
\noindent
{\bf Proof of Lemma \ref{lem401}.}  We can prove the case 
$ m=1 $ by a similar argument to the proof of Lemma \ref{lem309}.  
Next we consider the case $ m=2 $.  Put $ \uu_1^\e:=u_1^\e-M_1h^4 $.  
Since $ u_1^\e $ is a classical solution of (\ref{ellvi02}) with $ m=1 $ 
and satisfies $ u_1^\e\geq u_0^\e-M_1h^4 $ on $ \oom $, 
we observe that 
$$
   \frac{\uu_1^\e-u^\e_1}{h}-\uu^\e_{1,xx}+\be\uu_1^\e
    +\zeta_\e(\uu^\e_1-\varphi) 
    \leq \zeta_\e(u_1^\e-M_1h^4-\varphi) 
      -\zeta_\e(u_1^\e-\varphi)\leq  0\ \ {\rm in\ }\om 
$$
and $ \uu^\e_1(\pm 1)\leq u_1(\pm 1) $ for all $ \e\in (0,h^4) $ 
and $ h>0 $.  Applying the maximum principle, we have 
$ u_2^\e\geq u_1^\e-M_1h^4 $ on $ \oom $.  .  

By induction, we obtain the desired result.  $ \square $

\begin{lemma}
\label{lem405}
Assume $ q<r<qe $.  Then there exists $ h_3>0 $ such that 
$ u^\e_m>\varphi $ in $ (-1,0] $ and 
$ u_m^\e(0)\geq \varphi(0)+\sqrt{h}/4 $for all $ \e>0 $, 
$ m=1,2,\ldots,[T/h]$ and $ h\in (0,h_3) $.  
\end{lemma}

\noindent
{\bf Proof.}  
{\it Step 1.}  We claim that $ u^\e_m>\varphi $ in $ (-1,0] $ 
for all $ \e>0 $, $ m=1,2,\ldots,[T/h] $ and $ h>0 $.  

For each $ m=1,2,\ldots,[T/h] $, let $ \uU_m $ be the solution of 
$$
   \frac{\uU_m-\uU_{m-1}}{h}-\uU_{m,xx}+\be \uU_m=0\quad{\rm in\ }\om,\ 
    \uU_m(\pm 1)=\varphi(\pm 1).
$$
Since $ \uU_1 $ is a classical subsolution of 
(\ref{ellvi02}) with $ m=1 $, we have $ u^\e_1\geq \uU_1 $ 
on $ \oom $ by the maximum principle.  We see by induction that 
$ u^\e_m\geq \uU_m $ on $ \oom $ for all $ \e>0 $, 
$ m=1,2,\ldots,[T/h] $ and $ h>0 $.  Therefore, the claim of this step 
follows from $ \uU_m>0 $ in $ \om $ for all $ m=1,2,\ldots,[T/h] $ 
and $ \varphi\equiv 0 $ on $ [-1,0] $.  

{\it Step 2.} We prove that there exists $ h_4>0 $ such that 
$ u_m^\e(0)\geq \varphi(0)+\sqrt{h}/4 $
for all $ \e>0 $, $ m=1,2,\ldots,[T/h]$ and $ h\in (0,h_4) $.  

$ \uU_1 $ is given by 
$$
    \uU_1(x)=\frac{1}{h}\int_\om G_h(x,y)\varphi(y)dy
     +\frac{\varphi(1)\sh (z_h(x+1))}{\sh (2z_h)}, 
$$
where $ z_h:=\sqrt{\be+1/h} $, $ \sh(r):=\sinh(r) $ for $ r\in\bR $, 
\begin{equation}
\label{green}
   G_{a,h}(x,y):=\left\{\begin{array}{ll}
               \displaystyle{
     \frac{\sh(z_h(a-x))\sh(z_h(a+y))}{z_h\sh(2az_h)}}
          & (-a<y<x<a), \\
               \displaystyle{
     \frac{\sh(z_h(a+x))\sh(z_h(a-y))}{z_h\sh(2az_h)}}
          & (-a<x<y<a), 
              \end{array}\right. 
\end{equation}
for $ a>0 $ and $ G_h:=G_{1,h} $.  To estimate $ \uU_1(0) $, we directly 
calculate that  
for $ -1\leq x\leq 0 $, 
\begin{eqnarray}
\label{eq401}
   & & \uU_1(x)
    =A\sh(z_h(1+x))
      +\frac{h\sh (z_h(1+x))}{\sh (2z_h)}\left(
     \frac{qe^{\al+1}}{1+qh}-\frac{re^\al}{1+rh}\right),\\
\nonumber
   & & A:=\frac{1}{2hz_h\sh(2z_h)}
    \left\{\frac{e^{z_h}}{(z_h-\al-1)(z_h-\al)}
          -\frac{e^{-z_h}}{(z_h+\al+1)(z_h+\al)} \right\}.  
\end{eqnarray}
Putting $ x=0 $, we have 
$ u_1^\e(0)\geq \uU_1(0)\geq \sqrt{h}/3-C\sqrt{h}e^{-z_h} $ for all 
$ h>0 $.  In view of Lemma \ref{lem401} and $ \varphi(0)=0 $, 
selecting $ h_3>0 $ sufficiently small, we obtain the desired 
estimate.  $ \square $

\bigskip
The following lemma is suggested by the formal 
asymptotic expansion of solutions of (\ref{ellvi}) (cf. Section 6 below).  
\begin{lemma}
\label{lem402}
Put $ \mu_h:=\sqrt{h}-(\al+1/2) h$ and $ \rho:=(x-x_0-\mu_h)/2\sqrt{h} $.  
We define 
\begin{eqnarray}
\nonumber
  & &  \uu(x):=\left\{\begin{array}{l}
     \varphi(x)+e^{\al x}\{h^{3/2}w_3(\rho)+h^2(w_4(\rho)-M_5)\} 
       \quad  {\rm for\ }x\in [-1,x_0+\mu_h], \\
     \varphi(x)-M_5h^2e^{\al x}
      \hspace{46mm} {\rm for\ }x\in (x_0+\mu_h,1]
      \end{array}\right. \\
\label{func401}
   & & w_3(\rho):=r(e^{2\rho}-1-2\rho), \\
\nonumber
   & & w_4(\rho):=r\{e^{2\rho}-(1+2\rho+2\rho^2)
     +\al (e^{2\rho}-2\rho e^{2\rho}-1)\}.  
\end{eqnarray}
Then there are large $ M_5>0 $ and small $ h_4>0 $ such that 
$ \uu $ is a subsolution of {\rm (}\ref{ellvi02}{\rm )} with 
$ m=1 $ in the a.e.\;sense satisfying $ \uu(\pm 1)\leq \varphi(\pm 1) $ 
for all $ \e\in(0,h^4) $ and $ h\in (0,h_4) $.  
\end{lemma}

\noindent
{\bf Proof.}  Note that 
$ \uu\in W^{2,\infty}(\om)\cap C^2(\om\backslash\{x_0+\mu_h\}) $ 
and that $ w_3 $, $ w_4 $ satisfy 
\begin{eqnarray}
\label{ode401}
 & &  w_3-\frac{w_3^\pprime}{4}+r(2\rho+1)=0,\ 
   w_4-\frac{w_4^\pprime}{4}-\al w_3^\prime
         +r(2\rho^2+2\rho-\al)=0,\\
\label{esti421}
  & &  |\sqrt{h}w_3|+|hw_4|+|\sqrt{h}w_4^\prime|\leq C
   \quad{\rm in\ }(-(1+x_0+\mu_h)/2\sqrt{h},0].  
\end{eqnarray}

We divide our consideration into two cases.

{\it Case 1.}  $ x\in (-1,x_0+\mu_h] $. 

In this case, $ \rho\in (-(1+x_0+\mu_h)/2\sqrt{h},0] $.  
Using (\ref{ode401}), we compute that 
\begin{eqnarray*}
   \frac{\uu-\varphi}{h}-\uu_{xx}+\be\uu
    +\zeta_\e(\uu-\varphi) &\leq & e^{\al x}[h\{-M_5
    +r(\sqrt{h}w_3+hw_4)
    -\al \sqrt{h}w_4^\prime\} \\
   & & +qe^x-r-\sqrt{h}r(2\rho+1)-hr(2\rho^2+2\rho-\al)]. 
\end{eqnarray*}
We see from $ x_0=\log(r/q) $ and $ x-x_0=\mu_h+2\sqrt{h}\rho $ 
that 
\begin{eqnarray*}
   qe^x-r&=&r(e^{\mu_h+2\sqrt{h}\rho}-1)\leq r\left\{(\mu_h+2\sqrt{h}\rho)
  +\frac{1}{2!}(\mu_h+2\sqrt{h}\rho)^2+h^{3/2}\right\} \\
   &\leq & r\{\sqrt{h}(2\rho+1)+h(2\rho^2+2\rho-\al)\} \\
   & &   \quad+rh\left[\left(\al+\frac{1}{2}\right)
       \left\{-1-\sqrt{h}(2\rho+1)
    +\frac{h}{2}\left(\al+\frac{1}{2}\right)+\sqrt{h}\right\}\right].  
\end{eqnarray*}
Here we have used $ \mu_h+2\sqrt{h}\rho\leq\sqrt{h}/2 $ and 
the following inequality: 
$$
   e^\xi-1-\xi-\frac{1}{2}\xi^2\leq
     \left\{\begin{array}{ll}
      0 & {\rm if\ }\xi\leq 0,\\
      h^{3/2} & {\rm if\ }\xi\in [0,\sqrt{h}/2]
    \end{array}\right.
   \quad{\rm for\ small\ }h>0.  
$$
By the fact $ -\sqrt{h}(2\rho+1)\leq x_0 $, we get 
$$
   qe^x-r-r\{\sqrt{h}(2\rho+1)+hr(2\rho^2+2\rho-\al)\}\leq Ch.  
$$
From (\ref{esti421}) and this estimate, we have 
$$ 
  \frac{\uu-\varphi}{h}-\uu_{xx}+\be\uv+\zeta_\e(\uu-\varphi) 
     \leq he^{\al x}(-M_5+C)\quad{\rm in\ }(-1,x_0+\mu_h]
$$
for small $ \e>0 $ and $ h>0 $.  Taking $ M_5>0 $ large enough, 
we conclude that $ \uu $ is a classical subsolution of (\ref{ellvi02}) 
with $ m=1 $ in $ (-1,x_0+\mu_h) $ for small $ \e>0 $ and $ h>0 $.  

{\it Case 2.} $ x\in (x_0+\mu_h,1) $.

Taking $ M_5\geq 1$ and small $ h_5>0 $, we see that for all 
$ \e\in (0,h^4) $ and $ h\in (0,h_5) $, 
$$
  \frac{\uu-\varphi}{h}-\uu_{xx}+\be\uu+\zeta_\e(\uu-\varphi)
   \leq e^{\al x}(qe^x-r)-\zeta\left(-\frac{1}{h^2}\right)\leq 0
  \quad{\rm in\ }(x_0+\mu_h,1).  
$$

Therefore for large $ M_5>0 $ and small $ h_4>0 $, $ \uu $ is a 
subsolution of (\ref{ellvi02}) in the a.e.\;sense for 
$ \e\in (0,h^4) $ and $ h\in (0,h_4) $.  In view of (\ref{esti421}), 
we can get $ \uu(\pm 1)\leq \varphi(\pm 1) $ by replacing $ h_4 $ 
with a smaller one if necessary.  Thus the proof is completed.  $ \square $

\bigskip
\noindent
{\bf Proof of Theorem \ref{th405}.} It follows from Lemma \ref{lem402} 
and the maximum principle that $ u_1^\e\geq \uu $ on $ \oom $ for all 
$ \e\in(0,h^4) $ and $ h\in (0,h_5) $.  In view of Lemma \ref{lem405}, 
we have only to prove the assertion on $ [0,x_0+\sqrt{h}/2] $.  

First we treat the case $ m=1 $.  Let $ \rho $ and $ \mu_h $ be 
defined in Lemma \ref{lem402} and set $ \rho_1:=-(x_0+\mu_h)/2\sqrt{h} $, 
$ \rho_2:=-1/4+(2\al+1)\sqrt{h}/4 $.  
We observe by careful calculations that for small $ h>0 $, 
\begin{eqnarray*}
  & & \frac{d^2}{d\rho^2}
    \{h^{3/2}w_3+h^2(w_4-M_5)\}<0\quad{\rm on\ }[\rho_1,\rho_2],\\
  & & h^{3/2}w_3(\rho_1)+h^2(w_4(\rho_1)-M_5)\geq \frac{x_0}{4}h-M_5h^2,\\ 
  & & h^{3/2}w_3(\rho_2)+h^2(w_4(\rho_2)-M_5)\geq \frac{r}{10}h^{3/2}-M_5h^2.  
\end{eqnarray*}
Hence we have $ h^{3/2}w_3+h^2(w_4-M_5)\geq rh^{3/2}/20 $ on 
$ [\rho_1,\rho_2] $ and thus 
$$
   u_1^\e\geq\varphi+\frac{r}{20}h^{3/2}\quad 
    {\rm on\ }[0,x_0+\sqrt{h}/2]\ {\rm for\ any\ }\e\in(0,h^4)\ 
    {\rm and\ small\ } h>0.  
$$

In the case $ m\geq 2 $, Lemma \ref{lem401} and the above estimate yield 
that 
$$ 
   u_m^\e\geq\varphi+\frac{r}{20}h^{3/2}-M_1Th^3
   \geq\varphi+\frac{r}{40}h^{3/2}
   \quad{\rm on\ } [0,x_0+\sqrt{h}] 
$$ 
for all $ \e\in (0,h^4) $, $ m=1,2,\ldots,[T/h] $ and 
small $ h>0 $.  

Hence, selecting $ h_2>0 $ sufficiently small, 
we have $ u_m^\e\geq \varphi+ rh^{3/2}/40 $ on $ [0,x_0+\sqrt{h}/2] $ 
for all $ \e\in (0,h^4) $, $ m=1,2,\ldots,[T/h] $ and $ h\in (0,h_2) $.  
Thus we complete the proof.  
$ \square $

\bigskip
By Theorem \ref{th405}, 
we see that for any $ h\in (0,h_2) $, $ m=1,2,\ldots,[T/h] $ and 
$ \e\in (0,h^4) $,  
$ u_m^\e $ satisfies $ u_m^\e>\varphi $ in $ (-x_0,x_0) $ and thus 
$$
   \frac{u_m^\e-u_{m-1}^\e}{h}-u_{m,xx}^\e+\be u_m^\e=0\quad 
   {\rm in\ }(-x_0,x_0).  
$$
We prove Theorem \ref{th402}, based on this fact.  Before doing so, 
we give some preliminary analysis.   

In $ (-x_0,x_0) $, $ u_m^\e $ is given by 
$$
   u_m^\e(x)=\frac{1}{h}\int_\om G_{x_0,h}(x,y)u_{m-1}^\e(y)dy
     +\frac{u_m^\e(-x_0)\sh(z_h(x_0-x))}{\sh(2x_0z_h)}
     +\frac{u_m^\e(x_0)\sh(z_h(x_0+x))}{\sh(2x_0z_h)},  
$$
where $ G_{x_0,h} $ is defined by (\ref{green}) with $ a=x_0 $.  
In the sequel we set $ x_0=1 $ for simplicity.  Define 
$$
   \sG_h[\psi](x):=\frac{1}{h}\int_{-1}^1 G_h(x,y)\psi(y)dy
   \quad {\rm for\ }\psi\in C([-1,1]).  
$$
For $ -1\leq x\leq 0 $, we get from (\ref{eq401}) 
$$
  \sG_h[\varphi](x)
  =\frac{h\sh (z_h(1+x))}{\sh (2z_h)}\left(
     \frac{qe^{\al+1}}{1+qh}-\frac{re^\al}{1+rh}\right)
      -\frac{\varphi(1)\sh(z_h(1+x))}{\sh(2z_h)}+A\sh(z_h(1+x)).  
$$
On the other hand, we observe by tedious calculations 
that for $ 0<x\leq 1 $, 
\begin{eqnarray*}
 & & \sG_h[\varphi](x)
    =\varphi(x)+\frac{h\sh (z_h(x+1))}{\sh (2z_h)}\left(
     \frac{qe^{\al+1}}{1+qh}-\frac{re^\al}{1+rh}\right) \\
   & & \hspace{20mm} -he^{\al x}\left(\frac{qe^x}{1+qh}
    -\frac{r}{1+rh}\right)
    -\frac{\varphi(1)\sh(z_h(1+x))}{\sh(2z_h)}+B\sh(z_h(1-x)), \\ 
  & & B:=\frac{1}{2hz_h\sh(2z_h)}
    \left\{\frac{e^{z_h}}{(z_h+\al+1)(z_h+\al)}
          -\frac{e^{-z_h}}{(z_h-\al-1)(z_h-\al)} \right\}.  
\end{eqnarray*}
Noting that $ \varphi\equiv 0 $ on $ [-1,0] $ and that $ A>B>0 $, 
we have the following:  
\begin{eqnarray*}
  & & \sG_h[\varphi](x)\leq \varphi(x)+R(x)
    -\frac{\varphi(1)\sh(z_h(1+x))}{\sh(2z_h)}\quad
   {\rm for\ all\ }x\in[-1,1], \\ 
  & & R(x):=\frac{h\sh (z_h(1+x))}{\sh (2z_h)}\left(
     \frac{qe^{\al+1}}{1+qh}-\frac{re^\al}{1+rh}\right) 
   -\left(\frac{qhe^{(\al+1)x}}{1+qh}-\frac{rhe^{\al x}}{1+rh}\right)
     {\bf 1}_{\{x>0\}}\\
   & & \hspace{15mm}+\ds{A\sh(z_h(1-|x|))}.  
\end{eqnarray*}
Here $ {\bf 1}_{\{x>0\}}(x)=1 $ for $ x>0 $ and $ =0 $ for $ x\leq 0 $.  
Recalling $ u_0^\e=\varphi $, we see that 
\begin{eqnarray*}
   u_1^\e(x)&=&\sG_h[\varphi](x)
     +\frac{u_1^\e(-1)\sh(z_h(1-x))}{\sh(2z_h)}
     +\frac{u_1^\e(1)\sh(z_h(1+x))}{\sh(2z_h)}\\ 
    &\leq &u_0(x)+R(x) 
   +\frac{u_1^\e(-1)-u_0^\e(-1)}{\sh(2z_h)}\sh(z_h(1-x))
     +\frac{u_1^\e(1)-u_0^\e(1)}{\sh(2z_h)}\sh(z_h(1+x)).  
\end{eqnarray*}
We can inductively show that 
\begin{eqnarray*}
   u_m^\e(x)&\leq &u_{m-1}^\e(x)+\sG_h^{m-1}[R](x) 
   +\sum_{k=1}^m\frac{u_k^\e(-1)-u_{k-1}^\e(-1)}{\sh(2z_h)}
     \sG_h^{m-k}[\sh(z_h(1-\cdot))](x)\\
   & & +\sum_{k=1}^m\frac{u_k^\e(1)-u_{k-1}^\e(1)}{\sh(2z_h)}
     \sG_h^{m-k}[\sh(z_h(1+\cdot))](x), 
\end{eqnarray*}
where $ \sG_h^k[\psi]:=\sG_h[\sG_h^{k-1}[\psi]] $ and 
$ \sG_h^0[\psi]:=\psi $.  
Thus we need some pointwise estimates for $ \sG_{x_0,h}^m[R] $ and 
$ \sG_{x_0,h}^m[\sh(z_h(x_0\pm\cdot))] $ to prove Theorem \ref{th402}: 
\begin{theorem}
\label{th406}
For $ h>0 $, $ m=1,2,\ldots,[T/h] $ and $ x\in[-x_0,x_0] $, we have 
\begin{eqnarray}
\label{esti409} 
  & & \sG_{x_0,h}^m[\sh(z_h(x_0-|\cdot|))](x)
      \leq \sh(z_h(x_0-|x|))
        \sum_{k=0}^m\frac{a_{m,k}}{k!}(z_h|x|)^k+Ch,\\
\label{esti410}
  & &   \sG_{x_0,h}^m[\sh(z_h(x_0\pm\cdot))](x) 
    \leq \sh(z_h(x_0\pm x))
    \sum_{k=1}^m\frac{k}{2m-k}\frac{a_{m,k}}{k!}(z_h(x_0\mp x))^k, 
\end{eqnarray}
where $ a_{m,k}:=(2m-k)!/\{2^{2m-k}m!(m-k)!\} $.  
\end{theorem}
\begin{theorem}
\label{th407}
There exist $ x_4\in (0,x_0) $ and $ L_7 $, $ L_8 $, $ L_9>0 $ and 
$ h_5>0 $ such that for each $ h\in (0,h_5) $, 
\begin{eqnarray}
\label{esti404} 
  & & \sG_{x_0,h}^m[R](x)\leq \left\{\begin{array}{ll}
        \ds{L_7\sqrt{\frac{h}{m}}} & for\ x\in (-x_4,x_4),\\
        L_7h & for\ x=\pm x_4,
     \end{array}\right. 
   for\ all\ m=1,2,\ldots,[T/h], \\
\label{esti407}
  & & \sum_{k=1}^m
   \frac{\sG_{x_0,h}^{m-k}[\sh(z_h(x_0\pm\cdot))](x)}{\sh(2z_hx_0)} 
   \leq L_8h \\
\nonumber
  & & \hspace{30mm} 
   for\ all\ m=1,2,\ldots,[L_9/h|\log h|]\ 
    and \ x\in [-x_4,x_4].
\end{eqnarray}
\end{theorem}
We admit that Theorem \ref{th406} holds and prove 
Theorems \ref{th407} and \ref{th402}.  
We give the proof of Theorem \ref{th406} in Subsection 4.3 below.  

\bigskip
\noindent
{\bf Proof of Theorem \ref{th407}.}  
Set $ x_0=1 $ for notational simplicity.  Note that 
\begin{equation}
\label{ineq405}
     a_{m,k}\leq\frac{Ce^{-k^2/4(2m-k)}}{\sqrt{2m-k}}
   \quad{\rm for\ }k=0,1,\ldots,m.   
\end{equation}
This will be proved in Subsection 6.2 below.  

{\it Step 1.} We show (\ref{esti404}) for some $ x_{4,1}\in (0,1) $ 
and $ L_7>0 $.  

Since it is easily seen from (\ref{esti409}) that 
$$
   \sG_h^m[R](x)\leq \sqrt{h}e^{-z_h}\sh(z_h(1-|x|))
        \sum_{k=0}^m\frac{a_{m,k}}{k!}(z_h|x|)^k+Ch,  
$$
we have only to treat the first term of RHS of this inequality.  
Denote it by $ I_4(x) $.  

Using (\ref{ineq405}) and $ \ds{\sum_{k=0}^m\frac{(z_h|x|)^k}{k!}
\leq e^{z_h|x|}} $, we get  
\begin{eqnarray}
\label{esti415}
  I_4(x) & \leq & C\sqrt{\frac{h}{m}}\quad{\rm for\ all\ }
     m=1,2,\ldots,[T/h],\ x\in [-1,1]\ {\rm and\ }h>0.  
\end{eqnarray}

Fix $ x_{4,1}\in (0,1) $.  We consider only $ I_4(x_{4,1}) $ since 
$ I_4 $ is even.   We still denote it by $ I_4 $ if no confusion arises.  
Set $ m_1:=[2z_hx_{4,1}/5] $ and $ m_2:=3[z_hx_{4,1}] $.  Then we 
have from (\ref{ineq405}) 
$$
   I_4\leq C\sqrt{h}e^{-z_hx_{4,1}}
   \sum_{k=0}^{m}\frac{e^{-k^2/4(2m-k)}}{\sqrt{2m-k}}
    \frac{1}{k!}(z_hx_{4,1})^k.  
$$
We divide our consideration into three cases.  Let $ h_2 $ be given 
in Theorem \ref{th405}.  

{\it Case 1.} $ m\leq m_1 $.  

Using Stirling's formula, we get 
$$
   I_4\leq C\sqrt{h}e^{-z_hx_{4,1}}I_{4,1},\ 
    I_{4,1}:= \left\{1+\sum_{k=1}^{m}\frac{1}{\sqrt{k}}
     \left(\frac{z_hx_{4,1}e}{k}\right)^k\right\}.  
$$
Set $ \ga=k/z_hx_{4,1} $.  Then $ \ga\in (1/z_hx_{4,1},2/5) $ and 
$ (z_hx_{4,1}e/k)^k=\exp(z_hx_{4,1}(-\ga\log\ga+\ga)) $.  Since we see 
that $ -\ga\log\ga+\ga<4/5 $ for all $ \ga\in (1/z_hx_{4,1},2/5) $, 
we have 
$$
   I_{4,1}\leq Ce^{4z_hx_{4,1}/5}\left(1+\sum_{k=1}^m
     \frac{1}{\sqrt{k}}\right)
      \leq C(1+\sqrt{m_1})e^{4z_hx_{4,1}/5}\leq Ch^{-1/4}e^{4z_hx_{4,1}/5}.  
$$
Hence $ I_4\leq Ch^{1/4}e^{-z_hx_{4,1}/5} $ for all 
$ m=1,2,\ldots,[T/h] $ and $ h\in(0,h_2) $.  

{\it Case 2.} $ m>m_1 $.  

We may assume $ m>m_2 $.  Similar calculations as in Case 1 yield that 
\begin{eqnarray*}
   I_4 &\leq & C\sqrt{h}e^{-z_hx_{4,1}}\left(\sum_{k=0}^{m_1}
      +\sum_{k=m_1+1}^{m_2}+\sum_{k=m_2+1}^{m}\right)
     \frac{e^{-k^2/4(2m-k)}}{\sqrt{2m-k}}\frac{1}{k!}(z_hx_{4,1})^{k}\\
    &\leq & Ch^{1/4}e^{-z_hx_{4,1}/5}+C\sqrt{h}e^{-z_hx_{4,1}}
      \sum_{k=m_1+1}^{m_2}
     \frac{e^{-k^2/4(2m-k)}}{\sqrt{2m-k}}\frac{1}{k!}(z_hx_{4,1})^{k}.  
\end{eqnarray*}
Put $ \ds{I_{4,2}:=\sum_{k=m_1+1}^{m_2}\frac{e^{-k^2/4(2m-k)}}{\sqrt{2m-k}}
\frac{1}{k!}(z_hx_{4,1})^{k}} $.  

Setting $ m:=1/hs $, we observe that for $ k=m_1+1,\ldots,m_2 $, 
$$
   \frac{e^{-k^2/4(2m-k)}}{\sqrt{2m-k}}
   \leq \frac{C}{\sqrt{m}}e^{-m_1^2hs/8}\leq C\sqrt{hs}e^{-M_{4,1}s}
    \leq M_{4,2}\sqrt{h}.   
$$
Here and in the sequel $ M_{4,i} $'s ($ i\geq 1 $) are positive constants 
depending on $ x_{4,1} $.  Hence we have 
$$
   I_{4,2} \leq  M_{4,2}\sqrt{h}
      \sum_{k=m_1+1}^{m_2}\frac{(z_hx_{4,1})^{k}}{k!} 
    \leq M_{4,2}\sqrt{h}e^{z_hx_{4,1}}\quad{\rm for\ all\ }
     h\in (0,h_2).  
$$
Consequently we obtain $ I_4\leq M_{4,3}h $ for all $ m>m_1 $ and 
$ h\in (0,h_2) $.  

From Case 1 and 2, choosing $ L_7 $ large enough, we get 
(\ref{esti404}) for all $ h\in (0,h_2) $.  

{\it Step 3.} We show that for any $ x_{4,2}\in (0,3/4) $, there are 
$ h_{5,1}>0 $, $ L_8 $, $ L_9>0 $ such that (\ref{esti407}) holds 
for all $ m=1,2,\ldots,[L_9/h|\log h|] $, $ x\in [-x_{4,2},x_{4,2}] $ 
and $ h\in (0,h_{5,1}) $.  Denote by $ I_5(x) $ the left-hand side of 
(\ref{esti407}).  

Fix $ x_{4,2}\in (0,3/4) $.  For $ x\in [-x_{4,2},x_{4,2}] $, 
set $ m_1:=[2z_h(1-x)/5] $, $ m_2:=3[z_h(1-x)] $ and 
$ m_3:=[z_h^2(1-x)^2/100|\log h|] $.  We use (\ref{esti410}) and 
$ \sh(z_h(1+x))/\sh(2z_h)\leq 2e^{-z_h(1-x)} $ on $ [-1,1] $ to have 
\begin{eqnarray*}
   I_5(x) &\leq & Ce^{-z_h(1-x)}\sum_{k=1}^m\left[\sum_{l=1}^{m-k}
      \frac{le^{-l^2/4\{2(m-k)-l\}}}{\{2(m-k)-l\}^{3/2}}\right]
     \frac{(z_h(1-x))^l}{l!} \\
   &\leq & Ce^{-z_h(1-x)}\sum_{l=1}^m\left\{\sum_{k=l}^{m}\
     \frac{le^{-l^2/4(2k-l)}}{(2k-l)^{3/2}}\right\}  
     \frac{(z_h(1-x))^l}{l!}\quad 
    {\rm for\ all\ }x\in [-1,1].  
\end{eqnarray*}
We divide our considerations into two cases.

{\it Case 1.} $ m\leq m_1. $

It is easily observe from the fact $ le^{-l^2/4(2k-l)}/\sqrt{2k-l}\leq C $ 
for all $ k $, $ l\in\bN $ that 
$$
   \sum_{k=l}^{m}\frac{le^{-l^2/4(2k-l)}}{(2k-l)^{3/2}}
    \leq  C\sum_{k=l}^{m}\frac{1}{2k-l}
    \leq C\int_{1}^{T/h}\frac{1}{2r}dr \leq C|\log h| 
   \quad{\rm for\ all\ } l=1,2,\ldots,m.  
$$
Hence we use this inequality and the same argument as in Case 1 of 
Step 1 to obtain 
$$
    I_5(x)\leq C|\log h|e^{-z_h(1-x)}\sum_{l=1}^m\frac{(z_h(1-x))^l}{l!}
     \leq C|\log h|e^{-z_h(1-x)/5}
     \quad{\rm for\ all\ }x\in [-x_{4,2},x_{4,2}].  
$$

{\it Case 2.} $ m_1<m\leq m_3 $.  

We may consider $ m>m_2 $.  Similar calculations to Case 1 yield that 
\begin{eqnarray*}
   I_5(x) &\leq & Ce^{-z_h(1-x)}
       \left(\sum_{l=1}^{m_1}+\sum_{l=m_1+1}^{m_2}
        +\sum_{l=m_2+1}^{m}\right)
      \sum_{k=l}^{m}\frac{le^{-l^2/4(2k-l)}}{(2k-l)^{3/2}}
      \frac{(z_h(1-x))^l}{l!} \\
    &\leq & C|\log h|e^{-z_h(1-x)/5}+Ce^{-z_h(1-x)}
       \sum_{l=m_1+1}^{m_2}\sum_{k=l}^{m}\frac{le^{-l^2/4(2k-l)}}{(2k-l)^{3/2}}
      \frac{(z_h(1-x))^l}{l!}. 
\end{eqnarray*}
Set $ \ds{I_{5,1}(x):=\sum_{l=m_1+1}^{m_2}\sum_{k=l}^{m-1}
\frac{le^{-l^2/4(2k-l)}}{(2k-l)^{3/2}}\frac{(z_h(1-x))^l}{l!}} $.  
From the facts $ l\leq 2k-l\leq 2m_3 $ and $ l>m_1 $, we see that 
for $ k=l,\ldots,m $ and $ l=m_1+1,\ldots,m_2 $, 
$$
    \frac{le^{-l^2/4(2k-l)}}{(2k-l)^{3/2}}
      \leq \frac{1}{\sqrt{l}}e^{-m_1^2/8m_3}
    \leq Ch^{1/4}e^{-2|\log h|}=Ch^{9/4}. 
$$
Thus for $ x\in [-x_{4,2},x_{4,2}] $ and small $ h>0 $, 
\begin{eqnarray*}
   I_{5,1}(x)&\leq & C \sum_{l=m_1+1}^{m}\sum_{k=l}^m
    h^{9/4}\frac{(z_h(1-x))^l}{l!}
    \leq Ch^{9/4}m_3
     \sum_{l=m_1+1}^{m}\frac{(z_h(1-x))^l}{l!} \\
      &\leq & M_{5,1}h^{5/4}|\log h|e^{z_h(1-x)}\leq M_{5,1}he^{z_h(1-x)},  
\end{eqnarray*}
where $ M_{5,1} $ depends on $ 1-x_{4,2} $.  Consequently, we get 
$$
    I_{5,1}(x)\leq M_{5,1}h\quad{\rm for\ all\ }x\in [-x_{4,2},x_{4,2}]\ 
     {\rm and\ small\ }h>0.  
$$

Thus taking large $ L_8>0 $, 
$ L_9:=(1-x_{4,2})^2/100 $ and $ h_5:=\min\{h_2,h_{5,1}\} $, we obtain 
(\ref{esti407}) for all $ m=1,2,\ldots,[T/h] $, $ x\in [-x_4,x_4] $ 
and $ h\in (0,h_5) $.  

Setting $ x_4:=\min\{x_{4,1},x_{4,2}\} $, we complete the proof.  $ \square $

\bigskip
\noindent
{\bf Proof of Theorem \ref{th402}.}  {\it Step 1.}  
We claim that there are $ L_{4,1} $, $ L_{5,1}>0 $ and $ h_1>0 $ 
such that for all $ h\in (0,h_1) $, 
$ m=1,2,\ldots,[L_{5,1}/h|\log h|] $ and $ x\in\oom\backslash(-x_4,x_4) $.  
\begin{eqnarray}
\label{esti411}
   u_m^\e(x)-u_{m-1}^\e(x)\leq L_{4,1}h. 
\end{eqnarray}
Put $ h_1:=h_6 $.  It follows from Theorem \ref{th407} that 
\begin{equation}
\label{ineq410}
   u_m^\e(\pm x_4)-u_{m-1}^\e(\pm x_4)\leq (L_7+L_8)h
\end{equation}
for all $ m=1,2,\ldots,[L_9/h|\log h|] $ and $ h\in (0,h_1) $.  
Choosing $ L_4\geq L_7+L_8 $, we observe that $ u_0^\e+L_4h $ is a 
classical supersolution of (\ref{ellvi02}) in $ \om\backslash(-x_4,x_4) $ 
since $ u_0^\e(=\varphi) $ is smooth in this domain.  Hence we apply 
the maximum principle to have (\ref{esti411}) with $ m=1 $.  We 
inductively obtain(\ref{esti411}) for $ m=2,\ldots,[L_9/h|\log h|] $ 
and $ h\in (0,h_1) $.  Putting $ L_5:=L_9 $, we have the claim.  

{\it Step 2.} We derive the estimates of (1) and (2).  

From Theorem \ref{th407} and (\ref{esti411}), we obtain 
\begin{eqnarray}
\label{esti412}
  & & u_m^\e(x)-u_{m-1}^\e(x)\leq \left\{\begin{array}{ll}
   \ds{L_4\sqrt{\frac{h}{m}}} & {\rm if\ }|x|<x_4,\\
   L_4h & {\rm if\ }|x|\geq x_4 
      \end{array}\right. \\
\nonumber
  & & \hspace{20mm}{\rm for\ all\ }
   m=1,2,\ldots,[L_5/h|\log h|]\ {\rm and\ small\ }h\in (0,h_1).  
\end{eqnarray}
Hence (\ref{esti401}) holds for all $ m=1,2,\ldots,[L_5/h|\log h|] $.  

In the case $ m=[L_5/h|\log h|] $, we have 
$$
   u_{m}^\e-u_{m-1}^\e\leq L_4h\sqrt{\frac{|\log h|}{L_5}}
   \quad{\rm\ on\ } \oom.
$$
We can show by the maximum principle and induction that this estimate 
holds for $ m=[L_5/h|\log h|]+1,\ldots,[T/h] $.  Using (\ref{conv401}) 
and the above estimates, we have (\ref{esti401}) for all 
$ m=1,2,\ldots,[T/h] $ and $ h\in(0,h_1) $. 

The (\ref{esti402}) can be derived from (\ref{esti401}), (\ref{esti403}) 
and (\ref{esti405}).   The estimate of (2) is a consequence of 
(\ref{esti401}) and 
(\ref{conv401}).  
$ \square $

\subsection{Free boundary for the DMS-BS}

The problem (\ref{bs04}) has a unique free 
boundary.  However, it does not leads to the existence and uniqueness 
of that for the DMS-BS.  
To prove them is the purpose of this subsection.  

\begin{theorem}
\label{th404}
Assume $ q<r<qe $.  Then there is $ h_6>0 $ satisfying the following: 
For each $ h\in (0,h_6) $ and $ m=1,2,\ldots, [T/h] $, 
there exists a unique $ x_m\in (x_0,1] $ such that 
\begin{equation}
\label{free401}
  \{x\in\om\;|\; u_m(x)>\varphi(x)\}=(-1,x_m),\ 
  \{x\in\om\;|\; u_m(x)=\varphi(x)\}=[x_m,1].  
\end{equation}
Moreover, $ x_0<x_1\leq x_2\leq\cdots\leq x_m\leq\cdots\leq x_{[T/h]} $ 
and $ x_m\leq \min\{x_0+\sqrt{mh},1\} $ for all $ m=1,2,\ldots,[T/h] $.  
\end{theorem}

\noindent
{\bf Proof.}  Put $ h_6:=\min\{h_1,h_2,h_5\} $.  Notice by Theorems 
\ref{mono} and \ref{th405}, Lemma \ref{lem405} and (\ref{conv401}) 
that $ u_m>\varphi $ in $ (-1,x_0+\sqrt{h}/2) $ for all $ h\in (0,h_6) $.  
Since we easily observe by Theorem \ref{mono} that 
$ x_1\leq x_2\leq \cdots\leq x_m\leq\cdots $ if they exist, in the 
following we show the existence and uniqueness of $ x_m $.  

{\it Step 1.}  We treat the case $ m=1 $.  Set 
$ \rho=(x-x_0-\sqrt{h})/2\sqrt{h} $ and 
$$
    \ou_1(x):=\left\{\begin{array}{ll}
                \varphi(x)+e^{\al x}h^{3/2}w_3(\rho) 
             & {\rm for\ }x\in (0,x_0+\sqrt{h}), \\
                \varphi(x) & {\rm for\ }
                x\in [x_0+\sqrt{h},1),
            \end{array}\right.
$$
where $ w_3 $ is defined by (\ref{func401}).  We show that $ \ou_1 $ 
is a supersolution of (\ref{ellvi}) in (0,1) in the a.e.\;sense.  

Note that 
$ \ou_1\in W^{2,\infty}(\om)\cap C^2(\om\backslash\{x_0+\sqrt{h}\}) $ 
and that $ \ou_1\geq\varphi $ on $ [0,1] $.  We see by (\ref{ode401}) 
that in $ (0,x_0+\sqrt{h}) $,  
$$
  \frac{\ou_{1}-\varphi}{h}-\ou_{1,xx}+\be\ou_{1} 
    = e^{\al x}\{-h\al w_3^\prime+rh^{3/2}w_3
        +r(e^{\sqrt{h}(2\rho+1)}-1-\sqrt{h}(2\rho+1))\}.  
$$
It follows from the facts $ w_3 $, $ -w_3^\prime\geq 0 $ in 
$ (-\infty,0] $ and $ e^y\geq 1+y $ for all $ y\in\bR $ that 
$$
   \frac{\ou_{1}-\varphi}{h}-\ou_{1,xx}+\be\ou_{1} 
  \geq  0\quad {\rm in\ }(0,x_0+\sqrt{h}).  
$$
By (\ref{ineq304}) and this inequality, we see that $ \ou_1 $ is 
a supersolution of (\ref{ellvi}) in (0,1) in the a.e.\;sense.  

In view of $ \ou_1(0)\leq u_1(0) $, we modify $ \ou_1 $ to construct 
a viscosity supersolution of (\ref{ellvi}).  Put 
$ \eta:=2(x_0+\sqrt{T})^2 $.  Define $ \oW_1 $ by 
$$
   \oW_1(x)=\left\{\begin{array}{ll}
             -\ga(x-x_0)^3e^{-\eta/(x-x_0)^2} & (0\leq x\leq x_0), \\
              0 & (x_0<x\leq 1),  
          \end{array}\right.
$$
where $ \ga>0 $ is selected later.  Then $ \oW_1\in C^2(0,1) $, 
$ \oW_1\geq 0 $, $ \oW_{1,x}<0 $ in $ (0,1) $ and 
$ \oW_1(x_0)=\oW_{1,x}(x_0)=\oW_{1,xx}(x_0)=0 $.  
Moreover, it is easily observed by the choice of $ \eta $ that 
\begin{equation}
\label{ineq404}
   -\oW_{1,xx}+\be\oW_1\geq 0\quad{\rm in\ }(0,1).  
\end{equation}
Take $ \ga_1>0 $ satisfying $ \oW_1>\|\varphi\|_{L^\infty(\om)} $ 
on $ [0,x_0/4] $.  
Since $ \oW_2:=\|\varphi\|_{L^\infty(\om)} $ is a classical 
supersolution of (\ref{ellvi}) in $ \om $, setting 
$$
    \oU_1(x):=\left\{\begin{array}{ll}
           \oW_2 & {\rm if\ }-1\leq x\leq 0, \\
           \min\{\ou_1(x)+\oW_1(x),\oW_2\}& {\rm if\ }
                0<x\leq 1,
         \end{array}\right.
$$
we conclude that $ \oU_1 $ is a viscosity supersolution of 
(\ref{ellvi}) satisfying $ \oU_1(\pm 1)\geq\varphi(\pm 1) $.  

We have $ u_1\leq\oU_1 $ on $ [0,1] $ from the comparison principle 
for viscosity solutions.  Using this inequality, we can obtain a unique 
$ x_1 $ satisfying (\ref{free401}).  Indeed, $ \oU_1=\varphi $ on 
$ [x_0+\sqrt{h},1] $ 
implies that $ u_1=\varphi $ on $ [x_0+\sqrt{h},1] $.  Put 
$$ 
  x_1:=\inf\{y\in \om\;|\;u_1(x)=\varphi(x) 
  {\rm \ for\ all\ }x\in[y,1)\}(\leq x_0+\sqrt{h}).  
$$
Clearly $ u_1=\varphi$ on $ [x_1,1] $.  To verify $ u_1>\varphi $ in 
$ (x_0+\sqrt{h}/2,x_1) $, we suppose that there is 
$ \ox_1\in (x_0+\sqrt{h}/2,x_1) $ such that 
$ u_1(x)>\varphi(x) $ in $ (\ox_1,x_1) $ and $ u_1(\ox_1)=\varphi(\ox_1) $.  
Since $ u_1 $ and $ \varphi $ are solutions of (\ref{ellvi}) with $ m=1 $, 
we get $ u_1=\varphi $ on $ [\ox_1,x_1] $ by the uniqueness.  This 
contradicts to the definition of $ x_1 $ and hence $ u_1>\varphi $ 
in $ (x_0+\sqrt{h}/2,x_1) $.  This observation also leads to the 
uniqueness of $ x_1 $.  Therefore we have the desired result of Step 1.  

{\it Step 2.} We prove the case $ m=2 $.  

Let $ w_3 $ be defined by (\ref{func401}) and set 
$ \rho_m=(x-x_0-\sqrt{mh})/2\sqrt{mh} $ for $ m=1,2,\ldots,[T/h] $.  
By the facts $ w_3^\prime\leq 0 $ on  $ (-\infty, 0] $ 
and $ \rho_m\leq\rho_{m-1}\leq 0 $ on $ [0,x_0+\sqrt{(m-1)h}] $, 
we see that 
\begin{equation}
\label{ineq403}
    w_3(\rho_m)\geq w_3(\rho_{m-1})\quad 
  {\rm for\ all\ }x\in [-1,x_0+\sqrt{(m-1)h}]\ 
   {\rm and }\ m=1,2,\ldots,[T/h].  
\end{equation}
Define 
\begin{eqnarray*}
  & &  \ou_2(x):=\left\{\begin{array}{ll}
                \varphi(x)+e^{\al x}(2h)^{3/2}w_3(\rho_2) 
             & {\rm for\ }x\in [0,x_0+\sqrt{2h}], \\
                \varphi(x) & {\rm for\ }
                x\in [x_0+\sqrt{2h},1).
            \end{array}\right., \\
  & &  \oU_2(x):=\left\{\begin{array}{ll}
           \oW_2 & (-1\leq x\leq 0), \\
           \min\{\ou_2(x)+\oW_1(x),\oW_2\}& 
                (0<x\leq 1),
         \end{array}\right.
\end{eqnarray*}
We claim that $ \oU_2 $ is a viscosity supersolution of (\ref{ellvi}) 
with $ m=2 $.  

It is observed by (\ref{ineq404}) and 
$ u_1\leq \ou_1+\oW_1 $ in $ (0,x_0+\sqrt{h}] $ that 
\begin{eqnarray*}
  & & \frac{\oU_2-u_1}{h}-\oU_{2,xx}+\be\oU_2
      \geq\frac{\ou_2-u_1}{h}-\ou_{2,xx}+\be\ou_2 \\
  & &\hspace{43mm}\geq e^{\al x}\Big[\sqrt{h}
    \{(2\sqrt{2}-1)w_3(\rho_2)-w_3(\rho_1)\}
       -2h\al w_3^\prime(\rho_2) \\
  & & \hspace{49mm} 
       +r(2h)^{3/2}w_3(\rho_2)
        +r\{e^{\sqrt{2h}(2\rho_2+1)}-1-\sqrt{2h}(2\rho_2+1)\}\Big].  
\end{eqnarray*}
in $ (0,x_0+\sqrt{h}] $.  Using (\ref{ode401}), 
(\ref{ineq403}) with $ m=2 $ and $ e^y\geq 1+y $ for all 
$ y\in\bR $, we get 
$$
    \frac{\oU_2-u_1}{h}-\oU_{2,xx}+\be\oU_2 \geq 0 
  \quad{\rm on\ }(0,x_0+\sqrt{h}].  
$$
From $ \oU_2\geq \varphi $ on $ \oom $, (\ref{ineq304}) and this inequality 
it follows that $ \oU_2 $ is a viscosity supersolution of (\ref{ellvi}) 
with $ m=2 $.  

Thus we use the comparison principle for viscosity solutions to obtain 
$ u_2\leq \oU_2 $ on $ \oom $.  We can show by a similar argument to 
Step 1 that there exists a unique $ x_2 $ satisfying (\ref{free401}) 
with $ m=2 $ and $ x_2\leq x_0+\sqrt{2h} $.  

{\it Step 3.}  We consider the case $ m\geq 3 $.  

By induction we assume that there exists a unique $ x_{m-1} $ 
satisfying (\ref{free401}) and $ x_{m-1}\leq x_0+\sqrt{(m-1)h} $.  
Then define 
\begin{eqnarray*}
  & &  \ou_m(x):=\left\{\begin{array}{ll}
                \varphi(x)+e^{\al x}(mh)^{3/2}w_3(\rho_m) 
             & {\rm for\ }x\in (0,x_0+\sqrt{mh}], \\
                \varphi(x) & {\rm for\ }
                x\in [x_0+\sqrt{mh},1).
            \end{array}\right., \\
  & &  \oU_m(x):=\left\{\begin{array}{ll}
       \oW_2 & (-1\leq x\leq 0), \\
        \min\{\ou_m(x)+\oW_1,\oW_2\} & 
        (0<x\leq 1).
       \end{array}\right. 
\end{eqnarray*}
By a similar argument to Step 2, we can see that $ \oU_m $ is a 
viscosity supersolution of (\ref{ellvi}) and thus $ u_m\leq \oU_m $ 
on $ \oom $.  Therefore by the same way as in Step 1 we can find a 
unique $ x_m $ satisfying (\ref{free401}) and 
$ x_m\leq \min\{x_0+\sqrt{mh},1\} $.  $ \square $

\subsection{Proof of Theorem \ref{th406}.}

To prove Theorem \ref{th406}, we prepare some identities.  Define 
\begin{eqnarray*}
  & & \sI_k:= \int_{-1}^0(z_hy)^k\ch(2z_h(1+y))dy,\ 
  \sJ_k:=\int_0^x(z_hy)^k\ch(2z_hy)dy,\hspace{20mm}
\end{eqnarray*}
\begin{eqnarray*}
  & & \sK_k:= \int_x^1(z_hy)^k\ch(2z_h(1-y))dy,\ 
  \sL_k:=\int_{-1}^x(z_h(1-y))^k\ch(2z_h(1+y))dy,\\ 
  & & \sM_k:=\int_x^1(z_h(1-y))^k\ch(2z_hy)dy\quad 
   (k=1,2,\ldots,m).  
\end{eqnarray*}
Direct calculations yield that 
\begin{eqnarray*}
 & & \sI_0=\frac{1}{2z_h}\sh(2z_h),\ \sI_1=\frac{-1}{4z_h}(\ch(2z_h)-1),\\
  & & \sJ_0=\frac{1}{2z_h}\sh(2z_hx),\ 
    \sJ_1=\frac{1}{2z_h}\left\{P_{1,1}\sh(2z_hx)
     -P_{2,1}\ch(2z_hx)+\frac{1}{2}\right\},\\
  & & \sK_0=\frac{1}{2z_h}\sh(2z_h(1-x)),\ 
  \sK_1=\frac{1}{2z_h}\left\{P_{1,1}\sh(2z_h(1-x))+P_{2,1}
     \ch(2z_h(1-x))-P_{0,1}\right\},\\
   & & \sL_0=\frac{1}{2z_h}\sh(2z_h(1+x)),\
   \sL_1=\frac{1}{2z_h}\left\{Q_{1,1}\sh(2z_h(1+x))
    +Q_{2,1}\ch(2z_h(1+x))-Q_{0,1}\right\},\\
 & & \sM_0=\frac{1}{2z_h}(\sh(2z_h)-\sh(2z_hx)),\ 
  \sM_1=\frac{-1}{2z_h}\left\{Q_{1,1}\sh(2z_hx)
    +Q_{2,1}\ch(2z_hx)-Q_{0,1}\ch(2z_h)\right\},  
\end{eqnarray*}
where $ P_{0,1}:=1/2 $, $ P_{1,1}:=z_hx $, $ P_{2,1}:=1/2 $, 
$ Q_{0,1}:=1/2 $, $ Q_{1,1}:=z_h(1-x) $, $ Q_{2,1}:=1/2 $.  
For $ k\geq 2 $, the following identities hold.  
\begin{lemma}
\label{lem409}
Let $ \sI_k $, $ \sJ_k $, $ \sK_k $, $ \sL_k $, $ \sM_k $ be defined 
as above.  For $ k\geq 2 $, we have 
\begin{eqnarray*}
  & & (-1)^k\sI_k=\frac{k!}{2z_h}\left(-P_{0,k}
       +\frac{{\bf 1}_{k:{\rm even}}}{2^k}\sh(2z_h)
       +\frac{{\bf 1}_{k:{\rm odd}}}{2^k}\ch(2z_h)\right),\\  
  & & \sJ_k=\frac{k!}{2z_h}\left(P_{1,k}\sh(2z_hx)-P_{2,k}\ch(2z_hx)
      +\frac{{\bf 1}_{k:{\rm odd}}}{2^k}\right),\\
  & & \sK_k=\frac{k!}{2z_h}(P_{1,k}\sh(2z_h(1-x))
      +P_{2,k}\ch(2z_h(1-x))-P_{0,k}),\\  
  & & \sL_k=\frac{k!}{2z_h}(Q_{1,k}\sh(2z_h(1+x))
      +Q_{2,k}\ch(2z_h(1+x))-Q_{0,k}),\\
  & & \sM_k=\frac{-k!}{2z_h}\left\{Q_{1,k}\sh(2z_hx)
     +Q_{2,k}\ch(2z_hx)-\frac{{\bf 1}_{k:{\rm even}}}{2^k}
     \sh(2z_h)
     -\frac{{\bf 1}_{k:{\rm odd}}}{2^k}\ch(2z_h)\right\},  
\end{eqnarray*}
where 
\begin{eqnarray*}
  & & P_{0,k}:=\sum_{l=0\atop l:{\rm odd}}^k\frac{z_h^{k-l}}{2^l(k-l)!},\ 
   P_{1,k}:=\sum_{l=0\atop l:{\rm even}}^k \frac{(z_hx)^{k-l}}{2^l(k-l)!},\ 
   P_{2,k}:=\sum_{l=0\atop l:{\rm odd}}^k\frac{(z_hx)^{k-l}}{2^l(k-l)!},\\
  & & Q_{0,k}:=\sum_{l=0\atop{l:{\rm odd}}}^k\frac{(2z_h)^{k-l}}{2^l(k-l)!},\ 
 Q_{1,k}:=\sum_{l=0\atop{l:{\rm even}}}^k
     \frac{(z_h(1-x))^{k-l}}{2^l(k-l)!},\ 
   Q_{2,k}:=\sum_{l=0\atop{l:{\rm odd}}}^k
    \frac{(z_h(1-x))^{k-l}}{2^l(k-l)!}.  
\end{eqnarray*}
\end{lemma}

\noindent
{\bf Proof.} 
Integrating by parts we have, for $ k\geq 2 $, 
\begin{eqnarray*}
  & & (-1)^k\sI_k=\frac{1}{2z_h}\Bigg\{
     -\frac{1}{2}kz_h^{k-1}+\frac{z_h}{2}k(k-1)(-1)^{k-2}\sI_{k-2}\Bigg\},
    \hspace{50mm}
\end{eqnarray*}
\begin{eqnarray*}
  & & \sJ_k=\frac{1}{2z_h}\left\{(z_hx)^k\sh(2z_hx)
     -\frac{1}{2}k(z_hx)^{k-1}\ch(2z_hx)+\frac{z_h}{2}k(k-1)J_{k-2}
    \right\},\\  
  & & \sK_k=\frac{1}{2z_h}\Bigg\{(z_hx)^k\sh(2z_h(1-x))
    +\frac{1}{2}k(z_hx)^{k-1}\ch(2z_h(1-x))
      -\frac{1}{2}k(-z_h)^{k-1} \\
  & & \hspace{20mm} +\frac{z_h}{2}k(k-1)\sK_{k-2}\Bigg\}, 
\\
  & & \sL_k=\frac{1}{2z_h}\Bigg\{\{z_h(1-x)\}^k\sh(z_h(1+x))
     +\frac{k}{2z_h}\{z_h(1-x)\}^{k-1}\ch(z_h(1+x))
     -\frac{k}{2}(2z_h)^{k-1} \\
  & & \hspace{20mm} +\frac{z_h}{2}k(k-1)\sL_{k-2}\Bigg\},\\
  & & \sM_k=\frac{1}{2z_h}\Bigg\{-(z_h(1-x))^k\sh(2z_hx)
     -\frac{k}{2}(z_h(1-x))^{k-1}\ch(2z_hx)
     +\frac{z_h}{2}k(k-1)\sM_{k-2}\Bigg\}.  
\end{eqnarray*}
Using these recurrence formulae, we obtain the result.  $ \square $  

\bigskip
We separately prove (\ref{esti409}) and (\ref{esti410}) of Theorem 
\ref{th406}.  Put $ x_0=1 $ for the sake of simplicity.  

\bigskip
\noindent
{\bf Proof of (\ref{esti409}).}  
Set $ G_m(x):=\sG_h^m[\sh(z_h(1-|\cdot|))](x) $ for $ x\in[-1,1] $, 
$ \sh(r):=\sinh(r) $ and $ \ch(r):=\cosh(r) $ for $ r\in\bR $.  Note 
that 
\begin{eqnarray}
\label{eq432}
  & &  G_h>0\quad{\rm in\ }(-1,1)\times(-1,1),\ 
     \int_{-1}^1 G_h(x,y)dy\leq\frac{1}{hz_h^2}.  
\end{eqnarray}
In this proof we use the identities in Lemma \ref{lem409} and the 
following ones.  
\begin{eqnarray}
\label{eq425}
  & & \sh^2(z_h(1+y))=\frac{1}{2}(\ch(2z_h(1+y)-1),\\ 
\label{eq426}
  & & \sh(z_h(1-y))\sh(z_h(1+y))=\frac{1}{2}(\ch(2z_h)-\ch(2z_hy)),\\
\label{eq403}
  & & -\sh(z_h(1-x))\sh(2z_hx)+\sh(z_h(1+x))\sh(2z_h(1-x)) 
     =\sh(2z_h)\sh(z_h(1-x)), \\
\label{eq404}
  & & \sh(z_h(1-x))\ch(2z_hx)+\sh(z_h(1+x))\ch(2z_h(1-x)) 
     =\sh(2z_h)\ch(z_h(1-x)).  
\end{eqnarray}

{\it Step 1.} We estimate $ G_1(x) $.  

We calculate with using (\ref{eq425}) and (\ref{eq426}) 
to get for $ 0\leq x\leq 1 $, 
\begin{eqnarray*}
  G_1(x)& = & \frac{\sh(z_h(1-x))}{2hz_h\sh(2z_h)}\left\{
       \frac{1}{2z_h}\sh(2z_h)-1
        +x\ch(2z_h)-\frac{1}{2z_h}\sh(2z_hx)\right\} \\
  & & +\frac{\sh(z_h(1+x))}{2hz_h\sh(2z_h)}
      \frac{1}{2z_h}\sh(2z_h(1-x)).  
\end{eqnarray*}
By (\ref{eq403}), we get 
$$
   G_1(x)\leq\frac{b}{2hz_h^2}g_1(x)\quad{\rm for\ }x\in [0,1],
$$
where $ g_1(x):=(1+z_hx)\sh(z_h(1-x)) $ and $ b:=\ch(2z_h)/\sh(2z_h) $.  
By similar calculations we have $ G_1(x)\leq bg_1(-x)/2hz_h^2 $ for 
$ x\in [-1,0] $.  
Hence we obtain 
\begin{equation}
\label{ineq411}
    G_1(x)\leq \frac{b}{2hz_h^2}g_1(|x|)\quad{\rm for\ }x\in [-1,1].  
\end{equation}

{\it Step 2.} We consider the case $ m=2 $.  

It directly follows from (\ref{eq432}) and (\ref{ineq411}) that 
$ G_2(x)\leq b\sG_h[g_1](x)/2hz_h^2 $ for $ x\in [-1,1] $.  
We observe by (\ref{eq425}), (\ref{eq426}) and 
Lemma \ref{lem409} that for $ 0\leq x\leq 1 $, 
\begin{eqnarray*}
  \sG_h[g_1](x)& = &\frac{\sh(z_h(1-x))}{2hz_h\sh(2z_h)}
    \Bigg[\sI_0+(-1)\sI_1 
     +\frac{1}{z_h}\left\{z_hx+\frac{1}{2!}(z_hx)^2\right\}\ch(2z_h)
       -\{\sJ_0+\sJ_1\}\Bigg] \\
  & & +\frac{\sh(z_h(1+x))}{2hz_h\sh(2z_h)}
    \{\sK_0+\sK_1\} \\
  &=:&I_{5,1}.  
\end{eqnarray*}
We get from (\ref{eq403}) and (\ref{eq404}) 
\begin{eqnarray*}
   I_{5,1}&=&\frac{1}{2hz_h^2\sh(2z_h)}\left[
    \frac{1}{2}\{2+z_hx\}\sh(z_h(1-x))\sh(2z_h) \right. \\
   & &\hspace{25mm}\left. 
   +\left\{\frac{1}{4}+z_hx+\frac{1}{2!}(z_hx)^2\right\}
     \sh(z_h(1-x))(2\ch(2z_h)) \right.\\
   & & \hspace{25mm}\left. +\frac{1}{4}\left\{\ch(z_h(1-x))\sh(2z_h) 
     -\sh(z_h(1+x))\right\}\right].
\end{eqnarray*}
Since direct calculations yield that for all $ x\in [0,1] $, 
\begin{eqnarray}
\label{eq406}
  & & \ch(z_h(1-x))\sh(2z_h)-\sh(z_h(1+x))
    \leq \sh(z_h(1-x))\sh(2z_h)+\frac{1}{2}e^{-z_h},
\end{eqnarray}
we conclude that 
$$
  \sG_h[g_1](x)\leq \frac{b}{2hz_h^2}g_2(x)+\sP_2\quad 
   {\rm for\ }0\leq x\leq 1,  
$$
where 
$$
    g_2(x):=\frac{\sh(z_h(1-x))}{2hz_h^2}
    \left\{\frac{3}{2}+\frac{3}{2}z_hx+\frac{1}{2!}(z_hx)^2\right\},
   \ \sP_2:=\frac{e^{-z_h}}{8(2hz_h^2)\sh(2z_h)}.  
$$
We can obtain by the similar way as above 
$ \sG_h[g_1](x)\leq bg_2(-x)/2hz_h^2+\sP_2 $ for $ -1\leq x\leq 0 $.  
Consequently, we have 
\begin{eqnarray*}
  & & G_2(x)\leq\frac{b^2}{(2hz_h^2)^2}g_2(|x|)
    +b\sP_2\quad{\rm for\ }x\in [-1,1].
\end{eqnarray*}

{\it Step 3.} We estimate $ G_m(x) $ by induction.  

We assume that 
$$
   G_{m-1}(x)\leq  \frac{b^{m-1}}{(2hz_h^2)^{m-1}}g_{m-1}(|x|)
     +\sum_{l=2}^{m-1}b^{l-1}\sP_l\quad{\rm for\ }x\in [-1,1],  
$$
where $ g_{m-1} $ and $ \sP_l $ are defined by  
$$
 g_{m-1}(x):=\left(\sum_{l=0}^{m-1}\frac{c_{m-1,l}}{l!}(z_hx)^l\right)
        \sh(z_h(1-x)),\ 
  \sP_l:=\frac{e^{-z_h}}{4(2hz_h^2)^l\sh(2z_h)}
      \sum_{p=0}^{l-1} c_{l-1,p}
    \sum_{q=0\atop q:{\rm odd}}^p
     \frac{z_h^{p-q}}{2^q(p-q)!},
$$

We estimate $ \sG_h[g_{m-1}](x) $ for $ 0\leq x\leq 1 $.  
We calculate that 
\begin{eqnarray*}
  \sG_h[g_{m-1}](x)
      & \leq &\left[\frac{\sh(z_h(1-x))}{2hz_h\sh(2z_h)}
       \sum_{k=0}^{m-1}\frac{c_{m-1,k}}{k!}\left\{(-1)^k\sI_k
      +\frac{(z_hx)^{k+1}}{z_h(k+1)}\ch(2z_h)-\sJ_k\right\}\right. \\
\nonumber
   & & \left. +\frac{\sh(z_h(1+x))}{2hz_h\sh(2z_h)}
      \sum_{k=0}^{m-1}\frac{c_{m-1,k}}{k!}\sK_k\right]
      +\sum_{l=2}^kD^{l-1}\sP_l
    \quad{\rm for\ }0\leq x\leq 1.  
\end{eqnarray*}
It follows from Lemma \ref{lem409} that 
\begin{eqnarray*}
  I_{5,2}&:= & \sh(z_h(1-x))\{(-1)^k\sI_k-\sJ_k\}+\sh(z_h(1+x))\sK_k \\
  & =& \frac{k!}{2z_h}\Bigg[\sh(z_h(1-x))\left\{-P_{0,k}
       +\frac{{\bf 1}_{k:{\rm even}}}{2^k}\sh(2z_h)
       +\frac{{\bf 1}_{k:{\rm odd}}}{2^k}\ch(2z_h)\right. \\
  & & \hspace{20mm}  \left.-P_{1,k}\sh(2z_hx)
   +P_{2,k}\ch(2z_hx)-\frac{{\bf 1}_{k:{\rm odd}}}{2^k}\right\} \\
  & & \qquad  +\sh(z_h(1+x))\left\{P_{1,k}
    \sh(2z_h(1-x))+P_{2,k}\ch(2z_h(1-x))-P_{0,k}\right\}\Bigg],
\end{eqnarray*}
Using (\ref{eq403}), (\ref{eq404}) and $ \sh(2z_h)\leq \ch(2z_h)$, 
we obtain 
\begin{eqnarray*}
  I_{5,2}&\leq &\frac{k!}{2z_h}\Bigg[\sh(z_h(1-x))
      \left\{P_{1,k}\sh(2z_h)+\frac{1}{2^k}\ch(2z_h)
        \right\} 
    +P_{2,k}\ch(z_h(1-x))\sh(2z_h) \\
   & &   \qquad -\sh(z_h(1+x))P_{0,k}\Bigg]. 
\end{eqnarray*}
From (\ref{eq406}) and the fact $ P_{2,k}\leq P_{0,k} $ on 
$ [0,1] $ 
we get 
$$
   I_{5,2}\leq \frac{k!}{2z_h}\Bigg[\sh(z_h(1-x))
      \left\{(P_{1,k}+P_{2,k})\sh(2z_h)+\frac{1}{2^k}\ch(2z_h)\right\} 
    +e^{-z_h}P_{0,k}\Bigg].  
$$
Consequently, we have 
\begin{eqnarray*}
  \sG_h[g_{m-1}](x) &\leq &\frac{b\sh(z_h(1-x))}{2hz_h^2}
     \sum_{k=0}^{m-1}c_{m-1,k}
     \left\{\frac{1}{2}\left(\frac{1}{2^k}+\sum_{l=0}^k
     \frac{(z_hx)^{k-l}}{2^l(k-l)!}\right)
      +\frac{(z_hx)^{k+1}}{(k+1)!}\right\} \\
   & & +\frac{e^{-z_h}}{4\cdot 2hz_h^2\sh(2z_h)}
     \sum_{k=0}^{m-1}c_{m-1,k}\sum_{l=0\atop l:{\rm odd}}^k
       \frac{z_h^{k-l}}{2^l(k-l)!}+\sum_{l=2}^{m-1}b^{l-1}\sP_l.  
\end{eqnarray*}
Therefore setting 
\begin{eqnarray}
\label{eq407}
   g_m(x)&:=&\sh(z_h(1-x))\sum_{k=0}^m\frac{c_{m,k}}{k!}(z_hx)^k \\
\nonumber
     &=&\sh(z_h(1-x))\sum_{k=0}^{m-1}c_{m-1,k}
     \left\{\frac{1}{2}\left(\frac{1}{2^k}+\sum_{l=0}^k
     \frac{(z_hx)^{k-l}}{2^l(k-l)!}\right)
      +\frac{(z_hx)^{k+1}}{(k+1)!}\right\},  
\end{eqnarray}
we obtain 
$$
   G_m(x)\leq \frac{b^m}{(2hz_h^2)^m}g_m(x)
    +\sum_{k=2}^mb^{k-1}\sP_k\quad{\rm for\ }x\in [0,1].  
$$

Since we see 
$ \ds{G_m(x)\leq \frac{b^m}{(2hz_h^2)^m}g_m(-x)
    +\sum_{k=2}^mb^{k-1}\sP_k}
$ 
for $ -1\leq x\leq 0 $ by the same way as above, we conclude that 
$$
   G_m(x)\leq \frac{b^m}{(2hz_h^2)^m}g_m(|x|)
    +\sum_{k=2}^mb^{k-1}\sP_k\quad{\rm for\ }x\in [-1,1].  
$$

{\it Step 4.} We determine $ \{c_{m,k}\}_{k=0}^m $ for 
$ m=1,2,\ldots,[T/h] $.  

From (\ref{eq407}), we can obtain the following recurrence formulae:
for $ m=2,3,\ldots,[T/h] $ and $ k=2,3,\ldots,m-1 $, 
\begin{eqnarray}
\nonumber
  & & c_{1,1}=c_{1,0}=1,\ c_{m,m}=c_{m-1,m-1},\ 
   c_{m,m-1}=\frac{1}{2}c_{m-1,m-1}+c_{m-1,m-2}, \\
\label{eq421}
    & & c_{m,m-k}=\sum_{l=0}^{k}\frac{1}{2^{k-l}}c_{m-1,m-1-l},\ 
    c_{m,0}=\sum_{l=0}^{m-1}\frac{c_{m-1,l}}{2^l}.   
\end{eqnarray}
First, we easily get 
$$
    c_{m,m}=1,\ c_{m,m-1}=\frac{1}{2}(m+1)\quad{\rm for\ }m=1,2,\ldots,[T/h].  
$$
As for $ c_{m,m-2}$, using (\ref{eq421}) and these formulae, we have 
$$
   c_{m,m-2}=\frac{1}{2^2\cdot 2!}(m+1)(m+2).  
$$

We assume by induction that for $ m\geq 3 $, 
\begin{eqnarray*}
  & & c_{m-1,m-1-l}=\frac{1}{2^ll!}\prod_{p=1}^{l}(m-1+p)
   \quad {\rm for\ }l=1,2,\ldots,m-1.
\end{eqnarray*}
From (\ref{eq421}) and this equality we compute that for 
$ k=2,3,\ldots,m-1 $, 
\begin{eqnarray*}
   c_{m,m-k}&=&\frac{1}{2^k}\sum_{l=0}^{k}\frac{1}{l!}
      \prod_{p=1}^{l}(m-1+p)
     =\frac{1}{2^k}\left\{\frac{1}{2!}\prod_{p=1}^{2}(m+p)
      +\sum_{l=3}^{k}\frac{1}{l!}\prod_{p=1}^{l}(m-1+p)\right\} \\
     &=&\frac{1}{2^k}\left\{\frac{1}{3!}\prod_{p=1}^{3}(m+p)
      +\sum_{l=4}^{k}\frac{1}{l!}\prod_{p=1}^{l}(m-1+p) \right\}
    =\frac{1}{2^k k!}\prod_{p=1}^{k}(m+p).  
\end{eqnarray*}
Consequently, replacing $ k $ with $ m-k $, we obtain 
$$
   c_{m,k}=\frac{(2m-k)!}{2^{m-k}m!(m-k)!}
   \quad{\rm for\ }k=0,1,2\ldots,m. 
$$

{\it Step 5.}  We derive (\ref{esti409}).  

It is easy to see that for small $ h>0 $ and $ m=1,2,\ldots,[T/h] $, 
$ A\leq (1+Ce^{-3z_h})/z_he^{z_h} $ and $ (b/hz_h)^m \leq 1+e^{-2z_h} $. 
Besides, since $ c_{m,k}/2^m $ is the $ m $-th term of the binomial 
expansion of $ (1/2+1/2)^{2m-k} $, it is obvious that 
$ c_{m,k}/2^m\leq 1 $.  Using these facts, we get 
\begin{eqnarray*}
   \sG_h^m[\sh(z_h(1-|\cdot|))](x)&\leq & \frac{1}{z_he^{z_h}}g_m(|x|)
       +\frac{C}{z_he^{3z_h}}\sum_{k=0}^m\frac{z_h^k}{k!}
      +Ch+\sum_{l=2}^m\sP_l \\
    &\leq &\frac{1}{z_he^{z_h}}g_m(|x|)+Ch
       +\sum_{l=2}^m\sP_l  
\end{eqnarray*}
for small $ h>0 $.  Similarly we observe that 
$$
  \sum_{l=2}^m\sP_l \leq \frac{Ce^{-z_h}}{\sh(2z_h)}\sum_{l=2}^m
       \sum_{p=0}^{l-1}\frac{c_{l-1,p}}{2^{l+p}}
    \sum_{q=0}^p\frac{(2z_h)^q}{q!}
   \leq  Ce^{-z_h}\sum_{l=2}^m\sum_{p=0}^{l-1}\frac{1}{2^{p+1}}
  \leq Ch^{-1}e^{-z_h}
$$
for all $ m=1,2,\ldots,[T/h] $ and $ h>0 $.  Setting 
$ a_{m,k}=c_{m,k}/2^m $, we obtain (\ref{esti409}).  $ \square $

\bigskip
\noindent
{\bf Proof of (\ref{esti410}).}  
We treat only $ \sG_h^m[\sh(z_h(1+\cdot))] $ because 
$ \sG_h^m[\sh(z_h(1-\cdot))] $ can be similarly estimated.  
Set $ H_m(x):=\sG_h^m[\sh(z_h(1+\cdot))](x) $.  In this proof, we 
use the indentities in Lemma \ref{lem409}, (\ref{eq425}), (\ref{eq426}) 
and the following ones.  
\begin{eqnarray}
\label{eq429}
  & & \sh(z_h(1-x))\sh(2z_h(1+x))-\sh(z_h(1+x))(\sh(2z_h)-\sh(2z_hx))=0, \\
\label{eq431}
 & & \sh(z_h(1-x))(\ch(2z_h(1+x))-1)+\sh(z_h(1+x))(\ch(2z_hx)-\ch(2z_h))
  =0.  
\end{eqnarray}

{\it Step 1.}  We consider the case $ m=1 $.  

Using (\ref{eq425}) and (\ref{eq426}), we compute that 
\begin{eqnarray*}
 H_1(x) & =& \frac{\sh(z_h(1-x))}{2hz_h\sh(2z_h)}\cdot
   \frac{1}{2z_h}\sh(2z_h(1+x))\\
 & &  +\frac{\sh(z_h(1+x))}{2hz_h\sh(2z_h)}\left\{(1-x)\ch(2z_h)
    -\frac{1}{2z_h}(\sh(2z_h)-\sh(2z_hx))\right\}.  
\end{eqnarray*}
From (\ref{eq429}) we have 
\begin{equation}
\label{eq427}
 H_1(x)
   \leq  \frac{b}{2hz_h^2}h_1(x)\quad{\rm for\ }x\in[-1,1],\ 
    h_1(x):=z_h(1-x)\sh(z_h(1+x)).  
\end{equation}

{\it Step 2.} We estimate the case $ m=2 $.  

It follows from (\ref{eq432}) and (\ref{eq427}) that 
$ H_2(x)\leq b\sG_h[h_1](x)/2hz_h^2 $ for $ -1\leq x\leq 1 $.  
We see from (\ref{eq425}), (\ref{eq426}) and Lemma \ref{lem409} that 
\begin{eqnarray*}
 \sG_h[h_1](x) & \leq & 
  \frac{\sh(z_h(1-x))}{2hz_h\sh(2z_h)}
    \left[\frac{1}{2z_h}\left\{(z_h(1-x))\sh(2z_h(1+x))
    +\frac{1}{2}(\ch(2z_h(1+x))-1)\right\}\right] \\
   & &+\frac{\sh(z_h(1+x))}{2hz_h\sh(2z_h)}
     \left[\frac{1}{2!}z_h(1-x)^2\ch(2z_h)
      +\frac{1}{2z_h}\Bigg\{z_h(1-x)\sh(2z_hx)\right. \\
   & & \hspace{30mm}\left.
    +\frac{1}{2}(\ch(2z_hx)-\ch(2z_h))\Bigg\}\right].  
\end{eqnarray*}
We use (\ref{eq429}) and (\ref{eq431}) to obtain 
$$
  \sG_h[h_1](x)\leq  \frac{D}{2hz_h^2}h_2(x),\ 
   h_2(x):=\left\{\frac{1}{2!}(z_h(1-x))^2
    +\frac{1}{2\cdot 1!}(z_h(1-x))\right\}\sh(z_h(1+x)).  
$$
Consequenty we get 
$$
   H_2(x)\leq \frac{D^2}{(2hz_h^2)^2}h_2(x)\quad{\rm for\ }x\in[-1,1].  
$$

{\it Step 3.} We give an estimate for $ H_m(x) $ by induction.  

Suppose that for $ x\in [-1,1] $,
$$
  H_{m-1}(x)\leq \frac{D^{m-1}}{(2hz_h^2)^{m-1}}h_{m-1}(x),\ 
   h_{m-1}(x):=\sh(z_h(1+x))\sum_{p=1}^{m-1}
    \frac{d_{m-1,p}(z_h(1-x))^p}{2^{k-p} p!}, 
$$
It follows from (\ref{eq432}) and this inequailty that 
$ H_m(x)\leq b^{m-1}\sG_h[h_{m-1}](x)/(2hz_h^2)^{m-1} $.    
We easily see by (\ref{eq425}) and (\ref{eq426}) that 
\begin{eqnarray*}
  \sG_h[h_{m-1}](x)
  & \leq & \frac{\sh(z_h(1-x))}{2hz_h\sh(2z_h)}
     \sum_{k=1}^{m-1} \frac{d_{m-1,k}}{2^{m-1-k} k!}\sL_k\\
   & & +\frac{\sh(z_h(1+x))}{2hz_h\sh(2z_h)}
     \sum_{k=1}^{m-1} \frac{d_{m-1,k}}{2^{m-1-k} k!}
     \left\{\frac{z_h^k(1-x)^{k+1}}{k+1}\ch(2z_h)-\sM_k\right\}
\end{eqnarray*}
Using Lemma \ref{lem409}, (\ref{eq429}) and (\ref{eq431}), we have 
\begin{eqnarray*}
   & & \sh(z_h(1-x))\sL_k-\sh(z_h(1+x))\sM_k \\
  & & \quad =\frac{k!}{2z_h}\Bigg\{\sh(z_h(1+x))(Q_{1,k}\sh(2z_h)
     +Q_{2,k}\ch(2z_h))+(Q_{2,k}-Q_{0,k})\sh(z_h(1-x))\Bigg\} \\
  & & \quad =:I_6,  
\end{eqnarray*}
By 
$
   \ds{Q_{1,k}+Q_{2,k}=\sum_{l=0}^k\frac{\{z_h(1-x)\}^{k-l}}{2^l(k-l)!}} 
$ and $ Q_{0,k}\leq Q_{2,k} $, we get 
$$
  I_6\leq  \frac{bk!}{2z_h}\sh(2z_h)\sh(z_h(1+x))
    \sum_{l=0}^{k-1}\frac{(z_h(1-x))^{k-l}}{2^l(k-l)!}.
$$
Hence we obtain 
$$
  \sG_h[h_{m-1}](x)
   \leq  \frac{b\sh(z_h(1+x))}{2hz_h^2}\sum_{k=1}^{m-1} 
   \frac{d_{m-1,k}}{2^{m-k}}\left\{
     \sum_{l=0}^{k-1}\frac{(z_h(1-x))^{k-l}}{2^l(k-l)!}
    +\frac{2(z_h(1-x))^{k+1}}{(k+1)!}\right\}.  
$$

Therefore, setting 
\begin{eqnarray*}
   h_m(x)&:=&\sh(z_h(1+x))
      \sum_{k=1}^m\frac{d_{m,k}}{2^{m-k}k!}(z_h(1-x))^k \\
    &=& \sh(z_h(1+x))\sum_{k=1}^{m-1} \frac{d_{m-1,k}}{2^{m-k}} 
     \left\{
     \sum_{l=0}^{k-1}\frac{(z_h(1-x))^{k-l}}{2^l(k-l)!}
    +\frac{2(z_h(1-x))^{k+1}}{(k+1)!}\right\}, 
\end{eqnarray*}
we conclude that  
$$
   H_m(x)\leq \frac{b^m}{(2hz_h)^m}h_m(x)\quad{\rm for\ }
   x\in[-1,1]. 
$$

{\it Step 5.} We determine $ d_{m,k} $'s.  

It follows from the definition of $ h_m $ that for $ m=1,2,\ldots,[T/h] $ 
and $ k=2,3,\ldots,m-2 $, 
\begin{eqnarray}
\nonumber
  & & d_{m,m}=d_{m-1,m-1}=\cdots=d_{1,1}=1,\ d_{2,1}=1,\\ 
\label{eq430}
  & &  d_{m,m-1}
   =\sum_{l=1}^{m-1} d_{l,l},\ 
    d_{m,m-k}=\sum_{l=1}^{k+1}d_{m-1,m-l},
    \ d_{m,1}=\sum_{l=1}^{m-1}d_{m-1,l}.  
\end{eqnarray}
Here we see that 
$$
   d_{m,m-1}=m-1,\ 
   d_{m,m-2}=\frac{1}{2!}(m-2)(m+1).  
$$
We use these results to obtain $ d_{2,2}=1 $, $ d_{3,1}=d_{3,2}=2 $ 
and $ d_{3,3}=1 $.  

By induction we assume that for each $ m\geq 3 $, 
$$
   d_{m-1,m-1-k}=\frac{1}{k!}(m-1-k)\prod_{p=1}^{k-1}(m-1+p)\quad 
  {\rm for\ }k=2,3\ldots,m-2.  
$$
Then we calculate by using (\ref{eq430}) and this formula that for 
$ k=3,4,\ldots,m-2 $,
\begin{eqnarray*}
   d_{m,m-k}&=&1+(m-2)
    +\sum_{p=2}^{k}\frac{(m-1-p)}{(p-1)!}\prod_{q=0}^{p-2}(m+q)\\  
    &=&\frac{(m-2)(m+1)}{2!}
      +\sum_{p=3}^{k}\frac{(m-1-p)}{(p-1)!}\prod_{q=0}^{p-2}(m+q)\\ 
   &=&\frac{(m-3)}{3!}\prod_{p=1}^{2}(m+p)
      +\sum_{p=4}^{k}\frac{(m-1-p)}{(p-1)!}\prod_{q=0}^{p-2}(m+q)\\ 
   &=&\frac{(m-k)}{k!}\prod_{p=1}^{k-1}(m+p)=\frac{(m-k)(m+k-1)!}{m!k!}.  
\end{eqnarray*}

Replacing $ k $ with $ m-k $, we have 
$$
   d_{m,k}=\frac{k}{2m-k}\frac{(2m-k)!}{m!(m-k)!}\quad
   {\rm for\ }k=2,3,\ldots,m-2.  
$$
This formula clearly holds for $ k=1,m-1,m $.  Thus we obtain 
(\ref{esti410}).  $ \square $

\section{Proofs of main results}

First we prove Theorem \ref{th201}.  Let $ x_5:=\min\{x_1,x_4\} $ and 
$ \delta\in (0,T) $.  Set $ W_1:=[-x_5/2,x_5/2] $, $ W_2:=[-x_5,x_5] $, 
$ M_6:=\min\{L_5,x_5^2/16\} $ and 
$$
   P_h:=[0,M_6/|\log h|]\times W_1,\ 
   Q_{\delta,h}:=(\oom\times[0,T-\delta])\backslash
    ({\rm int\;}P_h\cup(\{0\}\times W_1)). 
$$
We show that there are $ K_1 $, $ K_2>0 $ and $ h_0>0 $ such that 
for all $ h\in (0,h_0) $,  
\begin{eqnarray}
\label{esti501}
  & &  \sup_{(t,x)\in P_h}|u(t,x)-u^h(t,x)|\leq K_1\sqrt{h},\\
\label{esti502}
  & &  \sup_{(t,x)\in Q_{\delta,h}}|u(t,x)-u^h(t,x)|\leq K_2\sqrt{h}|\log h|. 
\end{eqnarray}
Combining these estimates, we obtain the result of Theorem \ref{th201}.  

\bigskip
\noindent
{\bf Proof of (\ref{esti501}).}  
Choose $ h_{0,1}>0 $ so small that $ h<M_6/|\log h| $ for all 
$ h\in (0,h_{0,1}) $.  For $ t\in [0,h) $, Theorem \ref{th303} (2) 
directly yields that 
$$
  |u(t,x)-u^h(t,x)|=|u(t,x)-\varphi(x)|\leq C\sqrt{h}
   \quad{\rm for\ all\ }t\in [0,h),\ x\in\oom.  
$$
Hence in the following we consider the case 
$ t\in J_h:=[h,M_6/|\log h|] $.  The $ u $ and $ u^h $ are given by, 
respectively,  
\begin{eqnarray*}
  & & u(t,x)=[T(t)u_0](x)-\int_0^t E_y^{x_5}(t-s,x,x_5)u(s,x_5)ds 
   +\int_0^t E_y^{x_5}(t-s,x,-x_5)u(s,-x_5)ds,\\
  & & u^h(t,x)=\sG_{x_5,h}^m[u_0](x)+\sum_{k=1}^m\frac{u_k(-x_5)}{\sh(2x_5z_h)}
     \sG_{x_5,h}^{m-k}[\sh(z_h(x_5-\cdot))](x) \\
  & & \hspace{40mm}   +\sum_{k=1}^m\frac{u_k(x_5)}{\sh(2x_5z_h)}
     \sG_{x_5,h}^{m-k}[\sh(z_h(x_5+\cdot))](x),  
\end{eqnarray*}
for $ t>0 $, $ m=[t/h] $, $ x\in W_2 $ and $ h>0 $.  Here the family 
$ \{T(t)\}_{t\geq 0} $ is a contraction and analytic semigroup 
generated by the operator $ Au:=-u_{xx}+\be u $ in 
$ W_2 $ and $ D(A)=\{u\in C^2(W_2)\;|\;u(\pm x_1)=0 \} $ 
(cf. \cite[Corollary 3.1.21]{lu;95}).  We simply denote $ E^{x_5} $, 
$ \sG_{x_5,h} $ by $ \sG_h $, $ E $, respectively if no confusion 
arises.  

{\it Step 1.} We estimate $ \|[T(t)u_0]-\sG_h^m[u_0]\|_{C(W_2)} $.  

We use the contraction property of $ T(t) $ to have 
\begin{eqnarray*}
  \|[T(t)u_0]-\sG_h^m[u_0]\|_{C(W_2)}
   &\leq & \|[T(t-mh)u_0]-u_0\|_{C(W_2)} 
    +\|[T(h)[u_0]-\sG_h[u_0]\|_{C(W_2)} \\
   & & +\|[T((m-1)h)[\sG_h[u_0]]
        -\sG_h^{m-1}[\sG_h[u_0]]\|_{C(W_2)} \\
   &=:& I_{7,1}+I_{7,2}+I_{7,3}.  
\end{eqnarray*}
Since $ u_0(=\varphi) $ is Lipschitz on $ \oom $ and $ [T(t-mh)u_0] $ 
satisfies $ u_t-u_{xx}+\be u=0 $, it follows from the theory for 
parabolic equation that $ I_{7,1}\leq C\sqrt{h} $for all 
$ t\in [0,T) $, $ m=[t/h] $ and $ h>0 $.  In addition, direct 
calculations yield that for all $ h>0 $, 
$$
   I_{7,2}\leq \|[T(h)u_0]-u_0\|_{C(W_2)}
        +\|u_0(x)-\sG_h[u_0]\|_{C(W_2)}\leq C\sqrt{h}.  
$$
As for $ I_{7,3} $, we notice that $ \sG_h[u_0]\in D(A) $ and 
$ \|\sG_h[u_0]\|_{C^2(W_2)}\leq C/\sqrt{h} $.  Since it follows from 
\cite[Theorem 1.3]{ben;pau;04} that 
$ \|T((m-1)h)-\sG_h^{m-1}\|\leq Ch $, we get 
$$
    I_{7,3}\leq \|T((m-1)h)-\sG_h^{m-1}\|\|\sG_h[u_0]\|_{C^2(W_2)}
     \leq C\sqrt{h}.  
$$
for all $ m=1,2,\ldots,[T/h] $ and $ h>0 $.  Thus we obtain 
$$
   \sup_{t\in[h,T), m=[t/h]} 
    \|[T(t)u_0]-\sG_h^m[u_0]\|_{C(W_2)}\leq C\sqrt{h}\quad
   {\rm for\ all\ } h>0.  
$$

{\it Step 2.} We estimate 
$ \ds{I_{8,\pm}:=\int_0^t E_y(t-s,x,\pm x_5)u(s,\pm x_5)ds} $.  

We calculate that 
$$
   |E_y(t-s,x,\pm x_5)-E_{0,y}(t-s,x,\pm x_5)|
    \leq Ce^{-x_5^2/16(t-s)} 
$$
for all $ t,s\in (0,T) $ ($ t\ne s $) and $ x\in W_1 $.  
It is seen by this estimate that for $ t\in J_h $ and $ x\in W_1 $, 
\begin{eqnarray*}
   |I_{8,+}|&\leq &\|u\|_{C([0,T)\times\oom)}
     \int_0^t(|E_{0,y}(t-s,x,x_5)|+e^{-x_5^2/16(t-s)})ds \\
    & \leq & M_{7,1}\int_0^t \frac{x_5-x}{(t-s)^{3/2}}
       e^{-(x_5-x)^2/4(t-s)}ds +M_{7,2}h|\log h|^{-1}
     =: I_{8,1}+M_{7,2}h|\log h|^{-1}.  
\end{eqnarray*}
Here and in the sequel the constants $  M_{7,i} $'s ($ i\geq 1 $) 
depend on $ x_5 $, but not on $ h>0 $.  
Setting $ r=(x_5-x)/2\sqrt{t-s} $ and using 
$ \ds{\int_a^{+\infty}e^{-r^2}dr\leq e^{-a^2}/a} $ 
for $ a>0 $, we have 
$$
  I_{8,1}\leq M_{7,3}\int_{x_5/4\sqrt{t}}^{+\infty}e^{-r^2}dr 
     \leq  \frac{4M_{7,3}\sqrt{t}}{x_5}e^{-x_5^2/16t}\leq M_{7,4}h \quad
   {\rm for\ all\ }t\in J_h, x\in W_1\ {\rm and\ }h>0.
$$
Therefore we get $ I_{8,+}\leq M_{7,4}h $ for all $ t\in J_h $, 
We can get $ I_{8,-}\leq M_{7,5}h $ for all $ t\in J_h $, $ x\in W_1 $
and small $ h>0 $ by the same way as above.  Hence we have 
$$
   |I_{8,\pm}|\leq M_7h\quad {\rm for\ all\ }t\in J_h,\ x\in W_1\ 
  {\rm\ and\ small\ }h>0.  
$$
 
{\it Step 3.} We estimate $ \ds{I_{9,\pm}:=\sum_{k=1}^m
\frac{u_k(\pm x_1)}{\sh(2x_1z_h)}
     \sG_h^{m-k}[\sh(z_h(x_1\mp \cdot))](x)} $.  

It directly follows from the proof of Theorem \ref{th406} that 
$
   |I_{9,\pm}|\leq M_{7,6}h\ {\rm for\ all\ }m=1,2,\ldots,
    [M_6/h|\log h|],\ x\in W_1\ {\rm and}\ h>0.  
$

Therefore we have (\ref{esti501}) for $ h\in(0,h_{0,1}) $.  $ \square $

\bigskip
Next we prove (\ref{esti502}).  The point is to estimate the difference 
$ (u_m-u_{m-1})/h-u_t $.  To do so, we use the method similar to 
\cite{is;ko;91-2}, the precise comparison argument of viscosity 
solutions.  

Before proving (\ref{esti502}), we recall the definition and some 
elementary properties of the parabolic 2-jets.  Let $ W\subset\bR $ 
be an open interval.  For $ u:(0,T)\times W\to\bR $, we define 
$ \sP^{2,\pm}u(t,x) $ and $ \overline{\sP}^{2,\pm}u(t,x) $ as follows:  
\begin{eqnarray*}
  & & \sP^{2,+}u(t,x):=\left\{(a,p,X)\in\bR^3\;\left|
         \begin{array}{l}
           \ds{u(t+s,x+h)\leq u(t,x)+as+ph+\frac{1}{2}Xh^2}\\
           \hspace{30mm}+o(|s|+|h|^2)\quad
           {\rm as}\ (s,h)\to (0,0)
         \end{array}\right.\right\},\\
  & & \overline{\sP}^{2,+}u(t,x):=\left\{(a,p,X)\in\bR^3\;
      \left|\;\begin{array}{l}
            \exists \{(t_n,x_n,,a_n,p_n,X_n)\}_{n=1}^{+\infty}
      \subset(0,T)\times W\times\bR^3 
         \ {\rm such\ that} \\
            (t_n,x_n,u(t_n,x_n),a_n,p_n,X_n)\longrightarrow 
               (t,x,u(t,x),a,p,X)\\ 
          {\rm as\ }n\to+\infty\ {\rm and\ }
          (a_n,p_n,X_n)\in\sP^{2,+}u(t_n,x_n)
          \end{array}\right.\right\},\\
  & & \sP^{2,-}u(t,x):=\{(a,p,X)\in\bR^3\;|\;
          (-a,-p,-X)\in\sP^{2,+}(-u(t,x))\},\\
  & & \overline{\sP}^{2,-}u(t,x):=\{(a,p,X)\in\bR^3\;|\;
          (-a,-p,-X)\in\overline{\sP}^{2,+}(-u(t,x))\}.  
\end{eqnarray*}
We use the following lemma to obtain (\ref{esti502}).  
\begin{lemma}
\label{lem501}
Let $ u $, $ u_t\in C((0,T)\times W) $.  For any 
$ (t,x)\in (0,T)\times W $, if 
$ (a,p,X)\in\overline{\sP}^{2,+}u(t,x) $ {\rm (}or 
$ \overline{\sP}^{2,-}u(t,x) ${\rm )}, 
then $ a=u_t(t,x) $.  
\end{lemma}

\noindent
{\bf Proof.}  
Since $ u $ is differentiable with respect to $ t $, we can easily 
show that for any $ (t,x)\in (0,T)\times W $, if 
$ (a,p,X)\in \sP^{2,+}u(t,x) $, then $ a=u_t(t,x) $.  The assertion 
follows from the continuity of $ u_t $ and the definition of  
$ \overline{\sP}^{2,+}u(t,x) $.  

The case $ (a,p,X)\in\overline{\sP}^{2,-}u(t,x) $ is proved 
similarly.  $ \square $

\bigskip
\noindent
{\bf Proof of (\ref{esti502}).}  
First, we show that for any $ \delta>0 $, there exist $ K_{2,1}>0 $ 
and $ h_{0,2}>0 $ such that 
\begin{equation}
\label{esti505}
    \sup_{(t,x)\in Q_{\delta,h}}(u(t,x)-u^h(t,x))
      \leq K_{2,1}\sqrt{h}|\log h|\quad{\rm for\ all\ }h\in (0,h_{0,2}).
\end{equation}

{\it Step 1.} We define $ \uu^h(t,x) $ by 
$$
   \uu^h(t,x):=\left\{\begin{array}{ll}
       u_0(x) & {\rm for\ }t\in [0,h],\ x\in\oom,\\
       u_m(x) & {\rm for\ }t\in (mh,(m+1)h],\ m=1,2,\ldots,[T/h]
       \ {\rm and\ }x\in\oom.  
    \end{array}\right.
$$
For any $ \delta\in (0,T) $, put $ T_{\delta/2}:=T-\delta/2 $ and 
define 
$$
   \Phi(t,x,s,y):= u(t,x)-\uu^h(s,y)-\frac{1}{2\sqrt{h}}(t-s)^2 
     -\frac{1}{2\sqrt{h}}(x-y)^2
-\frac{\sqrt{h}}{T_{\delta/2}-t}-\frac{\sqrt{h}}{T_{\delta/2}-s}.
$$
Then $ \Phi $ is upper semicontinuous on 
$ Q_{\delta/2,h}\times Q_{\delta/2,h} $ and $ \Phi\longrightarrow -\infty $ 
as $ t $, $ s\nearrow T_{\delta/2} $.  Let $ (\ot,\ox,\os,\oy) $ 
be a maximum point of $ \Phi $ on $ Q_{\delta/2,h}\times Q_{\delta/2,h} $.  
We may consider $ \Phi(\ot,\ox,\os,\oy)\geq 0 $ 
because if otherwise, we easily get (\ref{esti505}) with 
$ K_{2,1}=4/\delta $.  

{\it Step 2.} We show that there is $ h_{0,3}>0 $ such that 
\begin{equation}
\label{esti515}
   u(\ot,\ox)-\uu^h(\os,\oy)\leq Ch^{1/4}|\log h|^{1/2}
   \quad{\rm for\ all\ }h\in (0,h_{0,3}).  
\end{equation}

For this purpose, we first study the behavior of $ (\ot,\ox,\os,\oy) $.  
It directly follows from $ \Phi(\ot,\ox,\os,\oy)\geq 0 $ that 
\begin{equation}
\label{ineq504}
   \frac{1}{2\sqrt{h}}(\ot-\os)^2
   +\frac{\sqrt{h}}{T_{\delta/2}-\ot}+\frac{\sqrt{h}}{T_{\delta/2}-\os}
   \leq u(\ot,\ox)-\uu^h(\os,\oy)\leq C.  
\end{equation}
Hence we get 
\begin{equation}
\label{esti509}
    |\ot-\os|\leq Ch^{1/4}.  
\end{equation}
Besides, since it is easily seen from 
$ \Phi(\ot,\oy,\os,\oy)\leq \Phi(\ot,\ox,\os,\oy) $ and Theorem 
\ref{th303} (2) that 
$$
   \frac{1}{2\sqrt{h}}(\ox-\oy)^2
     \leq u(\ot,\ox)-u(\ot,\oy)\leq C|\ox-\oy|,  
$$
we have 
\begin{equation}
\label{ineq301}
   |\ox-\oy|\leq C\sqrt{h}.  
\end{equation}

To obtain (\ref{esti515}), we divide our consideration into several 
cases.  Set $ \p_pQ_{\delta/2,h}:=\p Q_{\delta/2,h}
\backslash(\{T_{\delta/2}\}\times\oom) $.  

{\it Case 1.} $ (\ot,\ox) $ or $ (\os,\oy)\in \p_pQ_{\delta/2,h} $.  

We may assume $ (\ot,\ox)\in \p_pQ_{\delta/2,h} $ since the other case 
can be treated by the same way.  

{\it Subcase 1-1.} $ \ox\in\pom $ or ($ \ot=0 $ and $ |\ox|>x_1/2 $).  

Then $ u(\ot,\ox)=\varphi(\ox) $ and 
we get by Theorem \ref{mono} and (\ref{ineq301}) 
$$
   u(\ot,\ox)-\uu^h(\os,\oy)
   \leq \varphi(\ox)-\varphi(\oy) \leq C|\ox-\oy|\leq C\sqrt{h}.  
$$

{\it Subcase 1-2.} $ \ot\in[0,x_0^2/32|\log h|] $ and 
$ |\ox|\leq x_1/2 $.  

It follows from Theorem \ref{th402} (2), (\ref{esti501}), (\ref{esti509}) 
and (\ref{ineq301}) that for small $ h>0 $, 
\begin{eqnarray*}
  u(\ot,\ox)-\uu^h(\os,\oy)
  &\leq & u(\ot,\ox)-\uu^h(\ot,\ox)+\uu^h(\ot,\ox)-\uu^h(\os,\oy) \\
  &\leq & C\sqrt{h}+C(h^{1/4}\sqrt{|\log h|}+h^{1/2})
    \leq Ch^{1/4}\sqrt{|\log h|}.  
\end{eqnarray*}

{\it Case 2.} $ (\ot,\ox) $, $ (\os,\oy)\in {\rm int\;}Q_{\delta/2,h} $.  

Using the maximum principle for semicontinuous functions, we can 
find $ a $, $ b $, $ X $, $ Y\in\bR $ satisfying 
\begin{eqnarray}
\nonumber
  & & (a,(\ox-\oy)/\sqrt{h},X)\in\overline{\cal P}^{2,+}u(\ot,\ox),\ 
      (b,(\ox-\oy)/\sqrt{h},Y)\in\overline{\cal P}^{2,-}\uu^h(\os,\oy), 
    \\
\label{eq503}
  & & -\frac{3}{h}I\leq 
      \left(\begin{array}{cc}
        X & 0 \\
        0 & -Y
      \end{array}\right)\leq \frac{3}{h}
            \left(\begin{array}{cc}
        1 & -1 \\
        -1 & 1
      \end{array}\right),\ 
    a-b=\frac{\sqrt{h}}{(T_{\delta/2}-\ot)^2}
      +\frac{\sqrt{h}}{(T_{\delta/2}-\os)^2}.  
\end{eqnarray}
Set $ m=[\os/h] $.  Then $ \uu^h(\os,\oy)=u_m(\oy) $.  In addition, 
$ ((\ox-\oy)/\sqrt{h},Y)\in \overline{J}^{2,-}u_m(\oy) $ and 
$ X\leq Y $ (See \cite{cr;is;li;92} or \cite{ko;04} for the defintions 
of $ \overline{J}^{2,\pm} $ and the maximum principle for semicontinuous 
functions).  

We estimate the difference $ (u_m(\oy)-u_{m-1}(\oy))/h $.  It follows 
from $ \Phi_t(\ot,\ox,\os,\oy)=0 $, (\ref{eq503}), Theorem \ref{th305} 
and Lemma \ref{lem501} with $ W=\om $ that 
\begin{equation}
\label{eq501}
   a=u_t(\ot,\ox)=\frac{1}{\sqrt{h}}(\ot-\os)
     +\frac{\sqrt{h}}{(T_{\delta/2}-\ot)^2},\ 
    b=\frac{1}{\sqrt{h}}(\ot-\os)
     -\frac{\sqrt{h}}{(T_{\delta/2}-\os)^2}.  
\end{equation}
Noting $ b\geq 0 $ by Theorem \ref{mono}, we get $ \ot\geq \os $.  
We see from (\ref{esti314}), the first formula of (\ref{eq501}) and 
this fact that $ |\ot-\os|\leq C\sqrt{h|\log h|} $.  Substituting 
this into the second formula of (\ref{eq501}), we have 
$ 1/(T_{\delta/2}-\os)\leq Ch^{-1/4}\sqrt{|\log h|} $.  Thus 
$ \Phi(\ot,\ox,\os-h,\oy)\leq \Phi(\ot,\ox,\os,\oy) $, (\ref{eq501}) 
and this estimate yield that 
\begin{eqnarray}
\label{esti503}
    \frac{u_m(\oy)-u_{m-1}(\oy)}{h} \kern -1.5mm 
      &\leq & \kern -1.5mm \frac{1}{2h^{3/2}}\{(\ot-(\os-h))^2
      -(\ot-\os)^2\} \\
\nonumber
     & & 
    +\frac{h^{-1/2}}{T_{\delta/2}-(\os-h)}
      -\frac{h^{-1/2}}{T_{\delta/2}-\os} \\
\nonumber
      \kern -1.5mm & \leq &\kern -1.5mm 
     b+\frac{1}{2}\sqrt{h} 
          +\frac{h^{3/2}}{(T_{\delta/2}-\os)^3} 
    \leq b+C\sqrt{h}.  
\end{eqnarray}

By the way, since $ u $ is a viscosity subsolution of (\ref{bs04}) and 
$ u_m $ is a viscosity supersolution of (\ref{ellvi}), we have the 
following inequalities.  
\begin{eqnarray}
\label{ineq503}
  & & \min\{a-X+\be u(\ot,\ox),u(\ot,\ox)-\varphi(\ox)\}
    \leq 0, \\
\label{ineq501}
  & &   \min
   \left\{\frac{u_m(\oy)-u_{m-1}(\oy)}{h}-Y+\be u_m(\oy),
     u_m(\oy)-\varphi(\oy)\right\}\geq 0.  
\end{eqnarray}
If $ u(\ot,\ox)-\varphi(\ox)\leq 0 $ in (\ref{ineq503}), then 
we easily have by the above inequalities and (\ref{ineq301}) 
$$
   u(\ot,\ox)-u_m(\oy)\leq \varphi(\ox)-\varphi(\oy)\leq C\sqrt{h}.
$$
Thus, in the sequel we assume $ u(\ot,\ox)-\varphi(\ox)>0 $ for 
small $ h>0 $.  Then by (\ref{ineq503}), 
\begin{equation}
\label{ineq510}
   a-X+\be u(\ot,\ox)\leq 0\quad{\rm for\ small\ }h>0.  
\end{equation}
On the other hand, we easily get from (\ref{esti503}) and (\ref{ineq501}) that 
$$
      b+C\sqrt{h}-Y+\be u_m(\oy)\geq 0. 
$$
Combining (\ref{eq503}), (\ref{ineq510}) with this inequality, we have 
\begin{equation}
\label{esti516}
   u(\ot,\ox)-\uu^h(\os,\oy)=u(\ot,\ox)-u_m(\oy)\leq C\sqrt{h}.
\end{equation}

Taking $ h_{0,3}>0 $ small enough, we conclude that (\ref{esti515}) 
holds for all $ h\in (0,h_{0,3}) $.  

{\it Step 3.} We improve (\ref{esti515}) and establish (\ref{esti505}).  

Substituting (\ref{esti515}) into (\ref{ineq504}), we obtain 
\begin{equation}
\label{esti514}
  I_{10}:=\frac{1}{T_{\delta/2}-\ot}+\frac{1}{T_{\delta/2}-\os}
   \leq Ch^{-1/4}\sqrt{|\log h|}.  
\end{equation}
We observe from $ \Phi(\os,\ox,\os,\oy)\leq\Phi(\ot,\ox,\os,\oy) $, 
Theorem \ref{th303} (2) and this inequality that 
$$
    \frac{1}{2\sqrt{h}}(\ot-\os)^2\leq  C|\log h||\ot-\os|.  
$$
Hence we get $ |\ot-\os|\leq C\sqrt{h}|\log h| $.  By using this 
estimate, we improve the estimates in Case 1 as follows.  
\begin{equation}
\label{esti517}
   u(\ot,\ox)-\uu^h(\os,\oy)\leq C\sqrt{h}|\log h|^{3/2}\quad 
   {\rm if\ }(\ot,\ox)\ {\rm or\ }(\os,\oy)\in\p_pQ_{\delta/2,h}.  
\end{equation}
Therefore, (\ref{esti515}) can be improved in the following way:  
$$
   u(\ot,\ox)-\uu^h(\os,\oy)\leq C\sqrt{h}|\log h|^{3/2}.  
$$
Substituting this into (\ref{ineq504}) again, we get 
$ I_{10}\leq C|\log h|^{3/2} $.  Repeating the above argument, 
we have $ |\ot-\os|\leq C\sqrt{h|\log h|} $ and improve (\ref{esti517}) 
as  
$$
   u(\ot,\ox)-\uu^h(\os,\oy)\leq C\sqrt{h}|\log h|\quad 
   {\rm if\ }(\ot,\ox)\ {\rm or\ }(\os,\oy)\in\p_pQ_{\delta/2,h}.  
$$
Consequently, we have from (\ref{esti516}) and this estimate 
$$
   \Phi(\ot,\ox,\os,\oy)\leq C\sqrt{h}|\log h|\quad 
   {\rm for\ all\ }{\rm and\ }h\in (0,h_{0,3}).    
$$

Choosing a large $ K_{2,1}\geq C+4/\delta $, we obtain 
(\ref{esti505}).  

Next, we prove that for any $ \delta>0 $, there are $ K_{2,2}>0 $ and 
$ h_{0,4}>0 $ such that 
\begin{equation}
\label{esti508}
    \sup_{(t,x)\in Q_{\delta,h}}(u^h(t,x)-u(t,x))
 \leq K_{2,2}\sqrt{h}|\log h|\quad{\rm for\ all\ }h\in (0,h_{0,4}).  
\end{equation}

{\it Step 4.}  
Let $ u^h $ be defined by (\ref{approx01}).  For any 
$ \delta\in (0,T) $, define 
$$
   \Phi(t,x,s,y):= u^h(t,x)-u(s,y)-\frac{1}{2\sqrt{h}}(t-s)^2
     -\frac{1}{2\sqrt{h}}|x-y|^2
   -\frac{\sqrt{h}}{T_{\delta/2}-t}-\frac{\sqrt{h}}{T_{\delta/2}-s}. 
$$
Let $ (\ot,\ox,\os,\oy)\in Q_{\delta/2,h}\times Q_{\delta/2,h} $ 
be a maximum point of $ \Phi $.  We may consider 
$ \Phi(\ot,\ox,\os,\oy)\geq 0 $.  Note that (\ref{ineq504}), 
(\ref{esti509}) and (\ref{ineq301}) hold.  

If $ (\ot,\ox) $ or $ (\os,\oy)\in\p_pQ_{\delta/2,h} $, then we see by 
similar arguments to those in Case 1 of Step 2 that 
$$
    u^h(\ot,\ox)-u(\os,\oy)\leq Ch^{1/4}|\log h|\quad{\rm for\ small\ }h>0.  
$$
Thus we may assume $ (\ot,\ox) $, $ (\os,\oy)\in{\rm int\;}Q_{\delta/2,h} $ 
The maximum principle for semicontinuous functions 
yields $ a $, $ b $, $ X $, $ Y\in\bR $ satisfying 
\begin{eqnarray*}
  & & (a,(\ox-\oy)/\sqrt{h},X)\in\overline{\cal P}^{2,+}u^h(\ot,\ox),\ 
      (b,(\ox-\oy)/\sqrt{h},Y)\in\overline{\cal P}^{2,-}u(\os,\oy), 
    \\
  & & -\frac{3}{h}I\leq 
      \left(\begin{array}{cc}
        X & 0 \\
        0 & -Y
      \end{array}\right)\leq \frac{3}{h}
            \left(\begin{array}{cc}
        1 & -1 \\
        -1 & 1
      \end{array}\right),\ 
    a-b=\frac{\sqrt{h}}{(T_{\delta/2}-\ot)^2}
      +\frac{\sqrt{h}}{(T_{\delta/2}-\os)^2}.  
\end{eqnarray*}
Let $ m=[\ot/h] $.  Then note that $ u^h(\ot,\ox)=u_m(\ox) $ and hence 
$ ((\ox-\oy)/\sqrt{h},X)\in \overline{J}^{2,+}u_m(\ox) $ and $ X\leq Y $.  

We estimate $ (u_m(\ox)-u_{m-1}(\ox))/h $.  We see by 
$ \Phi_s(\ot,\ox,\os,\oy)=0 $ and Lemma \ref{lem501} that 
\begin{equation}
\label{eq502}
  a=\frac{1}{\sqrt{h}}(\ot-\os)+\frac{\sqrt{h}}{(T_{\delta/2}-\ot)^2},\   
  b=u_t(\os,\oy)=\frac{1}{\sqrt{h}}(\ot-\os)
     -\frac{\sqrt{h}}{(T_{\delta/2}-\os)^2}.
\end{equation}
Note $ \ot\geq\os $ by $ u_t(\os,\oy)\geq 0 $.  Dividing 
$ \Phi(\ot,\ox,\os,\oy)\geq \Phi(\ot-h,\ox,\os,\oy) $ by $ h $, 
we observe from this fact that 
\begin{eqnarray*}
\label{esti513}
   \frac{u_m(\ox)-u_{m-1}(\ox)}{h} 
\kern -1.5mm 
      &\geq & \kern -1.5mm \frac{1}{2h^{3/2}}\{(\ot-\os)^2
      -((\ot-h)-\os)^2\} \\
      & & 
    -\frac{h^{-1/2}}{T_{\delta/2}-(\ot-h)}
      +\frac{h^{-1/2}}{T_{\delta/2}-\ot} \\
      \kern -1.5mm & \geq &\kern -1.5mm 
     a-\frac{1}{2}\sqrt{h} 
          -\frac{h^{3/2}}{(T_{\delta/2}-\ot)^3}.  
\end{eqnarray*}
The (\ref{esti401}), (\ref{ineq504}) and the fact $ \ot\geq \os $ 
yield that $ 1/(T_{\delta/2}-\ot)\leq Ch^{-1/4}\sqrt{|\log h|} $ for 
small $ h>0 $.  Using this estimate, we have 
$$
    \frac{u_m(\ox)-u_{m-1}(\ox)}{h} \geq a-C\sqrt{h}\quad
   {\rm for\ small\ }h>0.  
$$
Since the remainder is totally similar to Step 1, we have 
$$
   u^h(\ot,\ox)-u(\os,\oy)\leq C\sqrt{h}|\log h|.  
$$
Thus taking $ h_{0,4}>0 $ small, we obtain (\ref{esti508}).  

Taking $ K_2:=\max\{K_{2,1},K_{2,2}\}+4/\delta $ and 
$ h_{0,2}:=\min\{h_{0,3},h_{0,4}\} $, we have the result.  
$ \square $

\bigskip
We establish the result of Theorem \ref{th201} by choosing 
$ h_0:=\min\{h_{0,1},h_{0,2}\} $.  

The proof of Theorem \ref{th202} is similar to \cite{jia;dai;04}, 
based on the limit operation of viscosity solutions due to 
\cite{ba;pe;87, ba;pe;88}.  

\bigskip
\noindent
{\bf Proof of Theorem \ref{th202}.}  Recall that 
$ x^h(t) $ is given by (\ref{approx01}).  Define 
$$
   \widehat x(t):=\limsup_{s\to t,h\to 0}x^h(s),\ 
   \ux(t):=\liminf_{s\to t,h\to 0}x^h(s).
$$
Notice $ \widehat x(0)=\ux(0)=x_0 $ by Theorem \ref{th404}. 

We show $ x^*\leq \ux $ in $ [0,T) $.  Fix $ t\in (0,T) $ and 
$ x>\ux(t) $.  Then there exist sequences $ \{h_n\}_{n=1}^{+\infty} $ 
and $ \{m_n\}_{n=1}^{+\infty} $ such that as $ n\to +\infty $, 
$$
  h_n\longrightarrow 0,\ m_nh_n\longrightarrow t,\ 
   x^{h_n}(s_n)=x_{m_n}\longrightarrow \ux(t).  
$$
Since $ x>x_{m_n} $ for large $ n\in\bN $, we get 
$ u^{h_n}(s_n,x)=u_{m_n}(x)=\varphi(x) $.  Letting $ n\to+\infty $, 
we have $ u(t,x)=\varphi(x) $ by Theorem \ref{th201} and thus 
$ x^*(t)\leq \ux(t) $.  

To prove $ \widehat x\leq x^* $ in $ [0,T) $, we suppose 
$ \widehat x(t_0)\geq x^*(t_0)+6\e_0 $ 
for some $ t_0\in (0,T) $, $ \e_0>0 $ and get a contradiction.  
By the continuity of $ x^* $ (cf. (\ref{free})), there exists 
$ \delta>0 $ such that 
\begin{equation}
\label{ineq505}
   \widehat x(t_0)>x^*(t)+5\e_0\quad{\rm for\ all\ }
   t\in (t_0-5\delta,t_0+5\delta).  
\end{equation}
Choose $ \{h_n\}_{n=1}^{+\infty} $ and $ \{m_n\}_{n=1}^{+\infty} $ 
satisfying 
\begin{equation}
\label{ineq502}
  h_n\longrightarrow 0,\ m_nh_n\longrightarrow t_0,\ 
   x^{h_n}(m_nh_n)\longrightarrow \widehat x(t_0)\quad
   {\rm as\ }n\to +\infty.  
\end{equation}
Take $ n_0\in\bN $ such that 
$$ 
  |m_nh_n-t_0|<\delta,\ 
  |x^{h_n}(m_nh_n)-\widehat x(t_0)|<\e_0\quad
    {\rm\ for\ all\ } n>n_0.  
$$
Using (\ref{free}), Theorem \ref{th404} and these facts, we observe 
that 
$$
   x^{h_n}(t)\geq x^{h_n}(m_nh_n)\geq x^*(t)+4\e_0\quad
   {\rm  for\ all\ }t\geq m_nh_n\ {\rm and\ } n>n_0.  
$$
This implies that for all $ n>n_0 $, 
\begin{equation}
\label{ineq511}
   u^{h_n}>\varphi\quad{\rm in\ }
   Q:=(t_0+2\delta,t_0+4\delta)\times(x^*(t_0)+\e_0,x^*(t_0)+3\e_0).  
\end{equation}
On the other hand, we notice $ u=\varphi $ in $ Q $.  

Fix $ y_0\in (x^*(t_0)+2\e_0,x^*(t_0)+4\e_0) $.  We derive 
\begin{equation}
\label{sub501}
    -\varphi_{xx}(y_0)+\be \varphi(y_0)\leq 0.   
\end{equation}
Let $ \phi=\phi(t,y) $ be a smooth function such that 
$ u-\phi $ takes its strict maximum at $ (t_0+2\delta/5,y_0) $ in $ Q $ 
and $ (\phi(t_0+3\delta,y_0),\phi_{xx}(t_0+3\delta,y_0)) 
=(\varphi(y_0),\varphi_{xx}(y_0)) $.  
Let $ (t_n,y_n) $ be a maximum point of $ u^{h_n}-\phi $ in $ Q $.  
Then it can be observed from Theorem \ref{th201} that 
\begin{equation}
\label{conv501}
  (t_n,y_n)\longrightarrow (t_0+3\delta,y_0),\ 
  u^{h_n}(t_n,y_n)\longrightarrow u(t_0+3\delta,y_0)=\varphi(y_0)
   \quad {\rm as}\ n\to \infty.  
\end{equation}
Put $ \widetilde m_n=[t_n/h_n] $.  Then from (\ref{ineq511}), 
$ u_{\widetilde m_n}(y_n)=u^{h_n}(t_n,y_n)>\varphi(y_n) $ for 
$ n>n_0 $.  Using the fact that $ u_{\widetilde m_n} $ 
is a viscosity subsolution of (\ref{ellvi}) with $ m=\widetilde m_n $ 
and Theorem \ref{mono}, 
we have the following inequality.
$$
    -\varphi_{xx}(y_n)
      +\be u_{\widetilde m_n}(y_n)\leq 0.  
$$
Letting $ n\to +\infty $, we get (\ref{sub501}).  

However, (\ref{sub501}) contradicts to (\ref{ineq304}) because 
of $ y_0>x^*(t_0)\geq x_0 $.  Thus we have $ \widehat x\leq x^* $ 
in $ [0,T) $ and conclude that $ \widehat x=\underline x=x^* $ in 
$ [0,T) $.  Applying \cite[Section 6]{cr;is;li;92}, we complete the 
proof.  $ \square $

\section{Appendix}
\subsection{Formal asymptotic expansion for (\ref{ellvi})}
This subsection is devoted to the formal asymptotic expansion of the 
solution of (\ref{ellvi}) with $ m=1 $ near the free boundary as 
$ h\searrow 0 $.  

Let $ x^*>0 $ be the free boundary of (\ref{ellvi}) with $ m=1 $.  
From the facts $ u_1=\varphi $ on $ [x^*,1] $ and $ u_1\in C^1(\om) $, 
it is sufficient to treat the following problem instead of (\ref{ellvi}):
\begin{equation}
\label{ode601}
     \ds{\frac{u-\varphi}{h}-u_{xx}+\be u=0} 
     \quad {\rm for\ }x<x^*,\ 
     u(x^*)=\varphi(x^*),\ u_x(x^*)=\varphi_x(x^*),\ 
\end{equation}
We rewrite (\ref{ode601}).  Set $ w(x):=(u(x)-\varphi(x))/e^{\al x} $.  
Then $ w $ satisfies
\begin{equation}
\label{ode602}
     \ds{-hw_{xx}-2\al hw_x+(1+rh)w=h(-qe^{x}+r)} 
     \quad {\rm for\ }x<x^*,\ 
     w(x^*)=w_x(x^*)=0.  
\end{equation}
The solution $ w $ of this problem is given by 
\begin{equation}
\label{sol601}
   w(x)=k_1e^{\lambda_+(x-x^*)}+k_2e^{\lambda_- (x-x^*)}
        +h\left(\frac{-qe^x}{1+qh}+\frac{r}{1+rh}\right)
   \quad{\rm for\ some\ }k_1,k_2\in\bR, 
\end{equation}
where $ \lambda_{\pm}:=-\al\pm z_h $.  Since we see from 
(\ref{ode602}) that $ w(x)=O(h) $ as $ h\to 0 $, we have $ k_2=0 $.  
Moreover, we observe by the conditions for $ u $ at 
$ x^* $ and Taylor expansion to $ \log(1+s) $ around $ s=0 $ 
that as $ h\searrow 0 $, 
\begin{eqnarray}
\label{esti601}
  x^*&=&x_0+\log\left\{\frac{(z_h-\al)(1+qh)}{(z_h-\al-1)(1+rh)}
        \right\}
    =x_0+\log\left(1+\frac{\sqrt{h}}{\sqrt{1+\be h}+\al\sqrt{h}}
      \right)\\
\nonumber
    &=&x_0+\frac{\sqrt{h}}{\sqrt{1+\be h}+\al\sqrt{h}}
      -\frac{1}{2}\left(\frac{\sqrt{h}}{\sqrt{1+\be h}+\al\sqrt{h}}
      \right)^2+O(h^{3/2})\\
\nonumber
   &=&x_0+\sqrt{h}-(\al+1/2)h+O(h^{3/2}),\\ 
\nonumber
\label{eq602}
  k_1&=&\frac{rh^{3/2}}{1+qh}+O(h^2). 
\end{eqnarray}
where $ x_0 $ is given by (\ref{free}).  Therefore using (\ref{sol601}) 
with these results, we have $ u(x)-\varphi(x)=O(h^{3/2}) $ near $ x<x^* $.  

To obtain the asymptotic expansion of $ w $ in terms of $ h $, from 
the above estimate, we may assume that $ w $ can be expanded as follows:
\begin{eqnarray}
\label{expa601}
  & & w(x)=h^{3/2}w_3(\rho)+h^2w_4(\rho)+O(h^{5/2}).  
\end{eqnarray}
Here $ \rho:=(x-x^*)/2\sqrt{h} $.  We impose the following from the 
conditions for $ w $ at $ x^* $ in (\ref{ode602}):
\begin{equation}
\label{cond601}
   w_i(0)=w_i^\prime(0)=0
   \quad {\rm for\ }i=3,4.  
\end{equation}
Substituting (\ref{expa601}) into (\ref{ode602}), we have by 
$ \be=\al^2+r $ and Taylor expansion to 
$ \exp(x_0+\sqrt{h}-(\al+1/2)h+O(h^{3/2})) $ as $ h \searrow 0 $ 
\begin{eqnarray}
\label{expa602}
   & & h^{3/2}\left(-\frac{w_3^\pprime}{4}+w_3\right) 
     +h^2\left(-\frac{w_4^\pprime}{4}-\al w_3^\prime
         +w_4\right) \\
\nonumber
   & & \qquad 
   =-r\{h^{3/2}(2\rho+1)+h^2(2\rho^2+2\rho-\al)\}
   +O(h^{5/2}).  
\end{eqnarray}
Here we have neglected $ k_1e^{\lambda_+y} $ since this is smaller 
than the last term of (\ref{sol601}) for $ x<x^* $ and $ h\searrow 0 $.  

We determine $ w_i $'s ($ i=3,4 $) from the above expansion.  
Comparing both sides of (\ref{expa602}), we can derive the following:
$$
   w_3-\frac{w_3^\pprime}{4}+r(2\rho+1)=0,\ 
   w_4-\frac{w_4^\pprime}{4}-\al w_3^\prime
         +r\left(2\rho^2+2\rho-\al\right)=0.  
$$
Solving these equations under (\ref{cond601}) we obtain 
\begin{eqnarray*}
  & & w_3(\rho)=r(e^{2\rho}-1-2\rho),\ 
  w_4(\rho)=r\{e^{2\rho}-(1+2\rho+2\rho^2)
    +\al(e^{2\rho}-2\rho e^{2\rho}-1)\}.  
\end{eqnarray*}

Therefore we conclude that as $ h\to 0 $, 
$$
    u(x)=\varphi(x)+e^{\al x}\left\{h^{3/2}
      w_3\left(\frac{x-x^*}{\sqrt{h}}\right)
     +h^2w_4\left(\frac{x-x^*}{\sqrt{h}}\right)+O(h^{5/2})\right\}
  \quad{\rm near\ }x<x^*.  
$$

\subsection{Proof of (\ref{ineq405})}

First, we may assume $ m\geq 3 $ and consider the case 
$ k=1,2,\ldots,m-1 $.  Because the case $ m=1,2 $ or $ k=0,m $ is easily 
proved.  We see from Stirling's formula that for all $ p\in\bN $, 
$$
    1+\frac{1}{12p}\leq \frac{p!}{\sqrt{2\pi}p^{p+1/2}e^{-p}}
    \leq 1+\frac{1}{12p}+\frac{C}{p^2}, 
$$
where $ C $ is independent of $ p $.  Using this inequality with 
$ p=m,2m-k,m-k $, we observe that for $ k=1,2,\ldots,m-1 $, 
\begin{equation}
\label{ineq601}
  \frac{(2m-k)!}{2^{2m-k}m!(m-k)!}
   \leq \frac{1+C/m}{\sqrt{(2m-k)\pi}}
     \left(\frac{m}{m-k/2}\right)^{-(m+1/2)}
     \left(\frac{m-k}{m-k/2}\right)^{-(m-k+1/2)}.  
\end{equation}

Set $ \delta=\delta(m,k):=k/\sqrt{2m-k} $.  Then we have 
$$
\left(\frac{m}{m-k/2}\right)^{m+1/2}
     \left(\frac{m-k}{m-k/2}\right)^{m-k+1/2} 
=\left(1+\frac{\delta}{\sqrt{2m-k}}\right)^{m+1/2}
\left(1-\frac{\delta}{\sqrt{2m-k}}\right)^{m-k+1/2}.  
$$
Denote the RHS of this formula by $ I_{11} $.  Since it follows from 
Taylor's expansion that 
$$
   \log(1+x)=x-\frac{1}{2}x^2+\frac{1}{3}x^3-\frac{x^4}{4(1+\theta)^4},\ 
    |\theta|<|x|<1, 
$$
we obtain for $ k=1,2,\ldots,m-1 $, 
\begin{eqnarray*}
 & & 
   \log I_{11}\geq\left(m+\frac{1}{2}\right)\left(\frac{\delta}{\sqrt{2m-k}}
      -\frac{\delta^2}{2(2m-k)}+\frac{1}{3}
   \frac{\delta^3}{(2m-k)^{3/2}}
        -\frac{1}{4}\frac{\delta^4}{(2m-k)^2}\right) \\
  & & \qquad +\left(m-k+\frac{1}{2}\right)\left(-\frac{\delta}{\sqrt{2m-k}}
      -\frac{\delta^2}{2(2m-k)}-\frac{1}{3}
   \frac{\delta^3}{(2m-k)^{3/2}} 
        -\frac{1}{4}\frac{\delta^4}{(2m-k)^2}\right) \\
  & & \quad=\frac{1}{2}\delta^2-\frac{\delta^2}{2(2m-k)}
    +\frac{\delta^4}{3(2m-k)}
    -\frac{\delta^4}{4(2m-k)}\left(1+\frac{1}{2m-k}\right). 
\end{eqnarray*}
Since it is easily seen by $ m\geq 3 $ that
$$
   \frac{1}{3}-\frac{1}{4}\left(1+\frac{1}{2m-k}\right)
    \geq \frac{1}{12}-\frac{1}{4m}\geq 0.  
$$
Thus we get 
$$
  \log I_{11}\geq \frac{1}{2}\delta^2-\frac{\delta^2}{2(2m-k)}
     >\frac{1}{4}\delta^2.  
$$
Therefore we obtain (\ref{ineq405}).  $ \square $

\bigskip
\noindent
{\bf Acknowledgement.} K. I. was partially supported by 
Grant-in-Aid for Scientific Research (No. 18204009, No. 20540117 and 
No. 20340026) of Japan Society of the Promotion of Science.  S. O. was 
partially supported by Grant-in-Aid for Scientific Research 
(No. 18340047 and No. 21654013) of Japan Society of the Promotion 
of Science.

\end{document}